\documentclass{article}
\usepackage{ifthen}
\usepackage{amsmath}
\usepackage{amssymb}
\usepackage{latexsym}
\usepackage{bbm}
\usepackage{pdfsync}
\usepackage{verbatim}
\usepackage{fancyhdr}
\usepackage{url}               
\usepackage[english]{babel}
\usepackage[applemac]{inputenc}
\usepackage{enumitem}
\usepackage{tikz}
\usepackage{pgfplots}
\usepackage{array,afterpage}   
\usepackage{array,booktabs,arydshln,xcolor}

\PassOptionsToPackage{dvipsnames,table}{xcolor}
 \usepackage{moreverb}
\usepackage{listings}
\usepackage{graphicx}
\usepackage{theorem}
\usepackage{mathtools}
\usepackage{mathabx}
\usepackage{amssymb}
\usepackage{amsmath}
\usepackage{geometry}
\usepackage[all]{xy}
\usepackage{graphicx}
\usepackage[pdftex,%
a4paper=true,%
bookmarks=true,%
bookmarksnumbered=true,%
colorlinks=true,%
urlcolor=blue,%
pdfstartview=FitH%
]{hyperref}

\newtheorem{theorem}{Theorem}[section]
\newtheorem{lemma}[theorem]{Lemma}

\newtheorem{proposition}[theorem]{Proposition}
\theoremstyle{definition}
\newtheorem{definition}[theorem]{Definition}
\newtheorem{example}[theorem]{Example}

\newtheorem{remark}[theorem]{Remark}

\numberwithin{equation}{section}
\numberwithin{figure}{section}

\newcommand{\bp}{{\it \text{Proof. }}}
\newcommand{\ep}{\hfill  $\square$}

\def\CC{\mathbb{C}}
\def\RR{\mathbb{R}}
\def\ZZ{\mathbb{Z}}
\def\NN{\mathbb{N}}
\def\eta{\mathcal{E}}
\def\mp{\mathcal{M}^{par}}
\def\smp{\mathcal{SM}^{par}}
\def\mps{\mathcal{M}_s^{par}}
\def\deg{\mathrm{deg}}
\def\pdeg{\mathrm{pardeg}}
\def\r{\mathrm{rank}}
\def\dim{\mathrm{dim}_{\CC}}

\def\en{\mathrm{End}}
\def\pe{\mathrm{parEnd}}
\def\spe{\mathrm{SparEnd}}
\def\l{\mathcal{L}}
\def\T{\mathcal{T}}
\def\M{\mathcal{M}_{\bullet}}
\def\H{\mathcal{H}}
\begin{document}
\hypersetup{
	colorlinks=true,     
	citecolor=red,
	linkcolor=blue     
	}
\thispagestyle{plain}
\fancyhead[L]{ }
\fancyhead[R]{}
\fancyfoot[C]{}
\fancyfoot[L]{ }
\fancyfoot[R]{}
\renewcommand{\headrulewidth}{0pt} 
\renewcommand{\footrulewidth}{0pt}
\newcommand{\montitre}{Parabolic Hitchin connection}
\newcommand{\auteur}{\textsc{ Zakaria Ouaras}} 
\newcommand{\affiliation}{Laboratoire J.-A. Dieudonn\'e \\
 Universit\'e C\^{o}te d'Azur\\
 Parc Valrose, 06108 Nice Cedex 02,   France\\
\url{ouaras@unice.fr}}

\begin{center}
{\bf  {\LARGE \montitre}}\\ \bigskip
{\large\auteur}\\ \smallskip \bigskip
 \affiliation \\ \bigskip
 \today \smallskip
\end{center}
 \begin{abstract}
In this paper,  we present an algebro-geometric construction of the Hitchin connection in the parabolic setting for a fixed determinant line bundle,  following \cite{pauly2023hitchin}.  Our strategy is based on Hecke modifications,  where we provide a decomposition formula for the parabolic determinant line bundle and the canonical line bundle of the moduli space of parabolic bundles.  As a special case, we construct a Hitchin connection on the moduli space of vector bundles with fixed non-trivial determinant.
 \end{abstract}
\pagestyle{fancy}
\fancyhead[R]{\thepage}
\fancyfoot[C]{}
\fancyfoot[L]{}
\fancyfoot[R]{}
\renewcommand{\headrulewidth}{0.2pt} 
\renewcommand{\footrulewidth}{0pt}
\section{Introduction}
Hitchin connection is a projective flat connection introduced by Hitchin  \cite{hitchin1990flat} in the context of geometric quantization of real symplectic manifolds,  motivated by the work of Witten \cite{witten1989quantum} in the study of quantum Chern-Simons theory.   For a given Lie group $G$, one associates the character variety of irreducible representation of the fundamental group of a topological surface $X$. This variety is canonically a symplectic manifold $(M,\omega)$ by the Atiyah-Bott-Goldman form \cite{Atiyah-Bott}, and a conformal structure on $X$ induces a complex structure on $M$.  Hitchin showed that the conditions for the existence of the projevtive connection are satisfied in this case for $G=\mathrm{SU}(r)$ and the closed oriented surface $X$ has genus $g \geq 2$ (except for $r$=$g$=$2$). 
By Narasimhan-Seshadri theorem, $M$ is identified with the moduli space of stable rank-$r$ vector bundles with trivial determinant over $X$  equipped with a complex structure,  thus a quasi-projective variety,  where the symplectic form $\omega$ is a K\"ahler form and the inverse of the determinant line bundle provides a pre-quantum line bundle (see \cite {MR783704} for details). \\

Scheinost-Schottenloher \cite{schottenloher1995metaplectic} generalize Hitchin's construction to higher dimensions,  on the relative moduli space of semistable holomorphic rank-r vector bundles $E$ with total Chern class one,  with trivial determinant over a family of K\"ahler varieties,  and also deals  with the case of parabolic vector bundles in dimension one.  This space is a compact complex variety equipped with a natural ample line bundle $\mathcal{L}$ that generalizes the determinant line bundle.  They show under the hypothesis that a square root of the canonical bundle $K_{\mathcal{M}/B}^{1/2}$ exists, that the pushforward sheaf $\mathcal{W}_k\footnote{The twist by $ K_{\mathcal{M}}^{1/2}$ is called the metaplectic correction}:= \varpi_*\left( \l^k \otimes K_{\mathcal{M}/B}^{1/2} \right)$ is locally free and equipped with projectively flat connection.  A special case is when $X$ is an \emph{elliptic surface},  Bauer \cite{bauer1989parabolic} give a description of the space of semistable rank-r bundles with total Chern class one and trivial determinant over $X$ in terms of parabolic bundles over the Riemann surface $C$ associated to $X$.
Hence, by this identification under the assumption, they get a projectively flat connection over the pushforward of the generalized determinant line bundle with a metaplectic correction in the parabolic setting \footnote{Note that we don't know how the generalized determinant bundle and the parabolic determinant bundle are related.}. Bjerre \cite{bjerre2018hitchin} removed the restriction requiring the existence of a square root of the canonical line bundle,  and proved the existence of the Hitchin connection over the space $\mathcal{W}^0_k:= \varpi_*\left( \l^k \right)$,  using a general construction of Hitchin connections in  geometric quantization,  as done in \cite{andersen2012hitchinT} and Hitchin connections in metaplectic quantization,  as done in \cite{andersen2012hitchin}.  \\

The constructions above were done using differential and K\"ahler geometric methods. Van Geemen-de Jong \cite{van1998hitchin} gave a purely algebraic approach for the construction of the Hitchin connection over a family $\mathcal{M}/S$ equipped with a line bundle $\mathcal{L}$. One of the conditions is the  equality $\mu_L \circ \rho=-\kappa_{\mathcal{M}/S},$
where $\kappa_{\mathcal{M}/S}$ is the Kodaira-Spencer map of a family $\mathcal{M}/S$,  $\mu_L$ is a map associated to a line bundle $L$ over $\mathcal{M}$ and $\rho$ is a symbol map.  We refer to several works related to algebro-geometric constructions of the Hitchin connection,  Faltings  \cite{faltings1993stable},  Axelrod-Witten-della-Pietra \cite{axelrod1991geometric}  and Ran \cite{ran2006jacobi}.    \\

Baier-Bolognesi-Martens-Pauly \cite{pauly2023hitchin} provide an algebo-geometric construction of the Hitchin connection over the relative moduli space of semi-stable rank-r vector bundles with trivial determinant over a family of complex projective curves of genus $g \geq 2$ (except $r$=$g$=2) parameterized by a variety $S$.  They use the so-called trace complex theory $\cite{beilinson1988determinant}$,  Bloch-Esnault quasi-isomorphism $\cite{bloch1999relative}$,  and Sun-Tsai isomorphism $\cite{sun2004hitchin}$.  They give a description of the Atiyah class of the determinant line bundle,  and the symbol map is a rational multiple of the quadratic part of the Hitchin system composed with the Kodaira-Spencer map of the family of curves.  
\paragraph{Setting of the problem:} Let $\pi_s: \mathcal{C} \longrightarrow S$ be a smooth versal family of projective curves of genus $g \geq 2$, parameterized by a projective variety $S$ and let $(\sigma_i)_{i\in I=\{1,2,...,N\}}$  $N$ disjoint sections of $\pi_s$,  i.e. 
$\forall i \neq j \in I$ and $\forall s \in S$,  we have   $\sigma_i(s) \neq \sigma_j(s)$, let denote $D=\sum_{i=1}^N \sigma_i(S)$ the associated divisor.  For a fixed  rank-$r$ parabolic type $\alpha_*=(k,\vec{a},\vec{m})$ with respect to the parabolic divisor $D$ and a relative line bundle $\delta \in \mathrm{Pic}^d(\mathcal{C}/S)$,  we denote by $\pi_e:\smp_{\mathcal{C}/S}:=\smp_{\mathcal{C}/S}(r,\alpha_*,\delta) \longrightarrow S$,  the relative moduli space 
of parabolic rank-$r$ vector bundles of determinant $\delta$ and type $\alpha_*$ equipped with the parabolic determinant line bundle$\Theta_{par} \in \mathrm{Pic}\left(\smp_{\mathcal{C}/S} /S \right)$.  We have the fibre product
\begin{equation}
 \def\cartesien{%
    \ar@{-}[]+R+<6pt,-1pt>;[]+RD+<6pt,-6pt>%
    \ar@{-}[]+D+<1pt,-6pt>;[]+RD+<6pt,-6pt>%
  }
\xymatrix{ 
\mathcal{C}\times_S\smp_{\mathcal{C}/S} \ar[rr]^{\pi_{n}} \ar[d]_{\pi_{w}} \cartesien &&  \smp_{\mathcal{C}/S} \ar[d]^{\pi_{e}} &  \\ 
\ \left( \mathcal{C},D \right) \ar[rr]_{\pi_{s}} && \ar@/^1pc/[ll]^{ \sigma_{i}} S
}
\label{fiber-product}
\end{equation}
\paragraph{Question:} Is there a projective flat connection on the vector bundle $\mathcal{V}^{par}(\alpha_*,\delta,\nu):=\pi_{e_*}\left( \Theta_{par}^{\nu} \right)$ over $S$ associated to a heat operator on the parabolic determinant line bundle ?
\paragraph{Motivation:} In the theory of conformal blocks Tsuchiya-Ueno-Yamada \cite{tsuchiya1989conformal} constructed for the Lie algebra $sl_{r}(\CC)$,  an integer $k$ and an N-tuple $\vec {\lambda}=\{\lambda_x\}_{x\in D}$ of dominant weights for $sl_{r}(\CC)$ over a marked curve $(C,D)$ a vector space  $\mathbb{V}_C(I,\vec{\lambda},k)$,  which can be glued to a vector  bundle $\mathbb{V}(I,\vec{\lambda},k)$ over the moduli space of marked curves $\mathfrak{M}_{g,N}$ parameterizing $N$-pointed curves of  genus-$g$,  called the vector bundle of conformal blocks.  In \cite{tsuchiya1989conformal} they constructed  a flat projective connection.  Moreover, Beauville-Laszlo  \cite{beauville1994conformal} proved that the vector spaces  $V_J$ in Hitchin's constructions are canonically identified with $\mathbb{V}_C(D,\vec{\lambda},k)$ over a curve with one point with trivial weights. It is natural to ask whether the connections of Hitchin  \cite{hitchin1990flat} and  Tsuchiya-Ueno-Yamada \cite{tsuchiya1989conformal}  coincide, this was proven by Laszlo \cite{laszlo1998hitchin}.  Pauly \cite{pauly1996espaces} gave a  generalization of Beauville-Laszlo's identification of the space of non-abelian theta functions $\mathrm{H}^0\left( \mathcal{SM}_C^{par},\Theta_{par}\right)$,  with the conformal blocks $\mathbb{V}(D,\vec{\lambda},k)$.  By \cite{tsuchiya1989conformal} there a projective flat connection on the conformal block side, hence it is natural to ask the above question, and to compare the two connections via Pauly's isomorphism.   \\

Biswas-Mukhopadhyay-Wentworth \cite{biswas2021ginzburg} gave a proof of the existence of the Hitchin connection for parabolic $G$-bundles following \cite{pauly2023hitchin} (we present here the case $G=SL_r$ \footnote{ Correspond to parabolic bundles with trivial determinant.}).  Their strategy uses Seshadri \cite{seshadri1977moduli}  identification of the moduli space of parabolic bundles with $\mathcal{SU}_{X/S }^{\Gamma}(r,d)$ the moduli space of $\Gamma$-bundle for a Galois group $\Gamma$ over a family of Galois cover $X/S$ of the family of curves $\mathcal{C}/S$  \footnote{See also \cite{mehta1980moduli},  \cite{bhosle1989parabolic}and \cite{BR93}.}.  Applying  \cite{pauly2023hitchin} Proposition 4.7.1 over the space $ \mathcal{SU}_{X/S} (r)$ applied to the determinant line bundle $\widehat{\l}$ and to conclude they proof that the map $\mu_{\widehat{\l}_Q}$ is an isomorphism hence the modified symbol map for the line bundle $(\widehat{\l}_Q)^{\nu}$ for a positive integer $\nu$ is given by 
\begin{equation}
\rho^{Hit}_{par,\Gamma}(\nu):=\mu^{-1}_{(\widehat{\l}_Q)^{\nu}} \circ \left(\cup[\widehat{\l}_Q] \circ \rho_{par} \circ \kappa_{\mathcal{C}/S} \right)
\label{BMWsymbol}
\end{equation}
where $\rho_{par}$ is the quadratic part of the parabolic Hitchin system and $\kappa_{\mathcal{C}/S}$ (resp.  $\kappa_{\smp_{\mathcal{C}}(r)/S}$) is the Kodaira-Spencer map of the family of curves (resp.  the Kodaira-Spencer map of the family of relative moduli space which depends on the Galois cover).  In a second paper $\cite{biswas2021geometrization}$, they prove that the symbol map $\rho^{Hit}_{par,\Gamma}(\nu)$  is independent on the parabolic weights in the case of full flag  using abelianization on parabolic Higgs bundles.  Hence, using the equality $\widehat{\l}_Q:=Q^*(\widehat{\l})\cong \Theta_{par}^{\vert \Gamma \vert /k},$ given in \cite{BR93}  (Proposition 4.14),  where $\vert \Gamma \vert$ is the order of $\Gamma$.  Set $\Theta$ the pull-back of the determinant line bundle by the forgetful map to $\mathcal{SU}_{\mathcal{C}/S}(r)$. Therefore, the symbol can be written as follows
$$\rho^{Hit}_{par,\Gamma}(\nu):=\vert \Gamma \vert \ \mu^{-1}_{\Theta_{par}^{\nu}} \circ \left(\cup[\Theta] \circ \rho_{par} \circ \kappa_{\mathcal{C}/S} \right).$$

Our work is independent of the work of \cite{biswas2021ginzburg}.  The object that we define are intrinsically attached to the marked curve  and the quasi-parabolic type.  Our proof is over $\smp_{\mathcal{C}/S}(r,\alpha_*,\delta)$,  the relative moduli space of parabolic rank-r vector bundles with fixed parabolic type $\alpha_*$, and fixed determinant $\delta \in \mathrm{Pic}^d(\mathcal{C}/S)$,  whereas \cite{biswas2021ginzburg} work for $\delta= \mathcal{O}_{\mathcal{C}}$.  We prove the following results:  \\

\noindent \textbf{Theorem \ref{Main Theorem 1}:} \emph{
Let $\nu$ be a positive integer. There exists a unique projective flat connection on the vector bundle ${\pi_e}_*(\Theta_{par}^{\nu})$ of non-abelian parabolic theta functions,  induced by a heat operator with symbol map $$ \rho^{Hit}_{par}(\nu):= \frac{1}{(\nu k+r)} \left( \rho_{par}\circ \kappa_{\mathcal{C}/S} \right).$$   }  

For $D=\emptyset$,  we have $\alpha_*$ = $k \in \NN^*$ the trivial parabolic type,  the space $\smp_{\mathcal{C}/S}(r,\alpha_*,\delta)$ coincides with $\mathcal{SU}_{\mathcal{C}/S}(r,\delta)$ the moduli space of semi-stable rank-r bundles with determinant $\delta$.  Hence,  $\Theta_{par}^{r/n}\equiv \l^k$ for $n:=\gcd(r,\deg(\delta))$,  where $\mathcal{L}$ is the determinant line bundle over $p_e:\mathcal{SU}_{\mathcal{C}/S}(r,\delta) \rightarrow S$, and $\rho_{par}=\rho^{Hit}$ is the Hitchin symbol map.  \\

\noindent \textbf{Theorem \ref{Main Theorem 2}:}  \emph{Let $ \mathcal{C}/S$ be a smooth versal family of complex projective curves of genus $g \geq 2$ (and $g \geq 3$ if $r=2$ and $\deg(\delta)$ is even).   
Let $k$ a positive integer.  There exists a unique projective flat connection on the vector bundle $p_{e_*}(\l^{k})$ of non-abelian theta functions,  induced by a heat operator with symbol map
$$ \rho(k):= \frac{n}{r(k+n)} \left( \rho^{Hit}\circ \kappa_{\mathcal{C}/S} \right).$$}

 Note that for $\delta$=$\mathcal{O}_{\mathcal{C}}$,  we recover the classical case proved in \cite{hitchin1990flat}.
\\
 
Our strategy is based on what we call Hecke maps,  which are rational maps given by associating to a parabolic bundle its Hecke modifications with respect to the elements of the filtrations,  hence we get vector bundles with fixed determinants.  Using this maps and the forgetful map over the moduli space of parabolic bundles,  we give in Proposition  \ref{Parabolic determinant bundle and Hecke modifications} a decomposition formula for the parabolic determinant line bundle $\Theta_{par}$
$$
\Theta_{par}^r:=  \lambda_{par}(\mathcal{E}_*)^r=\Theta^{n k} \otimes  \bigotimes\limits_{i=1}^N\bigotimes\limits_{j=1}^{\ell_i-1} \left( \Theta^n \otimes \Theta_{j}(i)^{n_j(i)} \right)^{p_j(i)}
$$
 and in the Proposition \ref{Canonical bundle} a decomposition formula for the relative canonical line bundle
$$K_{\smp_{\mathcal{C}/S}/S}=\Theta^{-2n} \otimes \bigotimes\limits_{i=1}^N \bigotimes\limits_{j=1}^{\ell_i-1}   \left( \Theta^n \otimes \Theta_{j}(i)^{n_j(i)} \right)$$
where $\Theta$ and $\Theta_j(i)$ are the pull-backs of the ordinary determinant line bundles under the forgetful map and  Hecke maps,  respectively.   We also prove a Beilinson-Schechtman type Theorem \footnote{See \cite{biswas2023BS} for a parabolic analog using another definition of parabolic Atiyah algebroid.}  \ref{Parabolic BSBE-2} over the moduli space of parabolic bundles,  where we express the Atiyah classes of the line bundles $\Theta$ and $\Theta_j(i)$ as a rational multiples of extension classes of exact sequences associated to a strongly-parabolic exact sequence (see definition \ref{SQPA}) of a parabolic bundle and its Hecke modifications.  We also define the parabolic Atiyah algebroid (see definition \ref{QPA}) and we show in Theorem \ref{infinitesimal deformations},  that the first cohomology group with values in this sheaf parameterizes the infinitesimal deformations of a marked curve equipped with a quasi-parabolic bundle.  \\  

Finally,  by proving the invariance of the parabolic Hitchin symbol $\rho_{par}$ under the Hecke modifications,  we obtain the van-Geemen-de Jong equation for any positive power of the parabolic determinant line bundle.  Additionally,  we establish the equation $\mu_{\Theta_{par}^{\nu}}= \left( \frac{\nu k+r}{k} \right) \cup [\Theta_{par} ]$ and we show the symbol map $\rho^{Hit}_{par}(\nu)$ is independent of the parabolic weights.   To conclude the proof,  we verify the remaining conditions of the existence and flatness of the connection using  Singh Theorem \ref{Singh} on Hitchin varieties.
\paragraph{Acknowledgements:}I extend my sincere gratitude to my thesis supervisor Christian Pauly,  for his invaluable guidance and constant encouragement throughout the completion of this work.  Additionally,  I would like to thank the reviewers of my thesis,  Xiaotao Sun and Richard Wentworth,  for their valuable feedback and suggestions.  Finally,  I am grateful to Michele Bolognesi and  Johan Martens for engaging discussions.
\section{Parabolic vector bundles and their moduli spaces}
Let $C$ be a smooth  projective complex  curve of genus $g \ge 2$ and $D$=$\{ x_1,x_2, ...,x_N \}$ a finite subset of points of $C$. 
The set $D$ will also be called a parabolic divisor. Set $I=\{1,2,...,N\}, $ where $N=\deg(D)$.
\subsection{Parabolic vector bundles}
\paragraph*{Parabolic type of a vector bundle} A parabolic type for a rank-$r$ vector bundle over the curve $C$ with respect to the parabolic divisor $D$ is the following numerical data $\alpha_*=(k,\vec{a},\vec{m})$ consisting of:  
\begin{itemize} 	
 	    \item A quasi-parabolic type  $\vec{m}=(\ell_i,m(i))_{i \in I}$,  where $\ell_i \in \NN^*$ are integers called the length,  and
  \begin{enumerate}
  	\item A sequence of integers called the flag type at $x_i\in D$
  	 $$m(i)=(m_1(i),m_2(i), ...,m_{\ell_ i}(i)),  \quad \ \mathrm{with} \quad m_j(i) \in \NN^*.  $$
  	\item We have for every $i \in I$  the relation: 
  	 $\sum\limits_{j=1}^{\ell_i} m_j(i)=r.$
    \end{enumerate} 
 \item A system of parabolic weights $(k,\vec{a})$,  where $k \in \NN^*$ and 
 $\vec{a}=(a_j(i))_{\substack{i \in I \\ 1 \leq j \leq \ell_i}} $
 a sequence of integers satisfying: 
 $$   0 \leq a_1(i) < a_2(i) <...<a_{\ell_i}(i)< k.$$ 
  \end{itemize} 
  We say that  $x_i \in D$ is a trivial point if $\ell_i=1$,  which implies that $m_1(i)=r$ and $m(i)=(r)$. \\
   We say that $\alpha_*$ is full flag parabolic type if $\ell_i=r$ for all $i\in I $,  thus $m_j(i)=1$  $\forall i,j$. 
  \begin{definition} [Parabolic vector bundles $\cite{seshadri1977moduli}$] Let $E$ be a rank-r vector bundle over $C$.  
A quasi-parabolic structure of quasi-parabolic  type $\vec{m}=(\ell_i,m(i))_{i\in I}$ on $E$ with respect to the parabolic divisor $D$, is given for each $i \in I$, by a linear  filtration of length $\ell_i$ on the fibre $E_{x_i}$
$$F^*_*(E): \ \ \ \ E_{x_i} = F^1_i(E) \supset F^2_i(E)  \supset \cdot\cdot\cdot\supset F^{\ell_i}_i(E) \supset F^{\ell_i+1}_i(E)=\{0\}
 $$
such that for $j \in \{1,2...,\ell_i\}$ we have 
  $\dim \left(F^{j}_i(E)/F^{j+1}_i(E)\right)$=$m_j(i)$.  A parabolic structure on $E$ with respect to the parabolic divisor $D$  is the data $(E,F^*_*(E),\alpha_*)$ where $\alpha_*=(k,\vec{a},\vec{m})$ is a fixed parabolic type,  and  $(E,F^*_*)$ is a quasi-parabolic structure over $E$ of type $\vec{m}$ with respect to the parabolic divisor $D$.  We denote a parabolic vector bundle by $E_*$. 
\end{definition}
For all $ i \in I$   and $ j \in \{1,2,...,\ell_i \}$,  we define the following quotients $Gr^j_i(E):=\left( F^{j}_{i}(E)/F^{j+1}_{i}(E) \right)$
and $Q^j_i(E):=\left( E_{x_i}/F^{j+1}_{i}(E) \right)$, with dimensions $m_j(i)$ and $r_j(i)=\sum\limits_{q=1}^{j} m_q(i)$,  respectively. 
\begin{definition}[parabolic degree, slope and Euler characteristic]  \item Let $E_*$ be a parabolic bundle over $C$.  We define:
\begin{enumerate}
\item The parabolic degree:  $\pdeg(E)=\deg(E)+\frac{1}{k} \sum\limits_{i=1}^N  \sum\limits_{j=1}^{\ell_i} m_j(i)  a_j(i).$
 \item The parabolic slope: $ \mu_{par}(E)=\pdeg(E)/\r(E).$
 \item The parabolic Euler characteristic: $ \chi_{par}(E)=  \pdeg(E)+\r(E)(1-g).  $
\end{enumerate}
 \end{definition}	    
 \begin{definition}[parabolic and strongly parabolic endomorphisms] Let $E_*$ be a parabolic bundle over $C$  and let $f \in \en(E)$. Then, $f$ is a  parabolic  (res.  strongly parabolic) endomorphism, if for all $i \in I$ and $ j \in \{1,2,...,\ell_i \}$ one has
 $$  f_{x_i} \left(F^{j}_i(E)\right) \subset F^{j}_i(E) \quad res.  \quad f_{x_i} \left(F^{j}_i(E)\right) \subset F^{j+1}_i(E) .$$
We denote, the sheaf of parabolic endomorphisms and  strongly parabolic endomorphisms by $\pe(E)$ and $\spe(E)$,respectively, which are locally free. By definition, we have the inclusions
\begin{equation} 
 \spe(E)  \hookrightarrow \pe(E) \hookrightarrow \en(E).
 \label{SPE}
\end{equation}
 \end{definition}
 \begin{proposition} $\cite{yokogawa1991moduli}$ Let $E_*$ be a parabolic vector bundle over $C$ with respect to the parabolic divisor $D$. Then, we have a canonical isomorphism of locally free sheaves 
 $$ \pe(E)^{\vee} \cong \spe(E) \otimes \mathcal{O}_C(D).$$
 This isomorphism is given by the non-degenerate trace paring:  
 $$\begin{array}{ccccc}
 \mathrm{Tr} & : & \pe(E)\otimes \spe(E)& \longrightarrow &  \mathcal{O}_C(-D) \\
 & &\phi \otimes \psi & \longmapsto &  \mathrm{Tr}(\phi\circ \psi). 
\end{array}$$
 \end{proposition}
by dualizing \eqref{SPE} we get
 $$ \en(E)^{\vee}\cong \en(E)  \hookrightarrow  \spe(E)(D)  \hookrightarrow  \pe(E)(D) .$$
 \subsection{Moduli spaces of parabolic bundles} 
 To construct the moduli space of parabolic vector bundles over a curve $C$, we need a notion of semi-stability and stability which will depend on the parabolic type $\alpha_*$. Mehta-Seshadri constructed the moduli space of semistable parabolic vector bundles over a smooth family  of projective complex curve $C$. In this subsection,  we recall the existence theorem of a coarse moduli space parameterizing parabolic bundles. The main reference is  $\cite{seshadri1982fibers}$.
\begin{theorem} [Mehta-Seshadri  \cite{mehta1980moduli}] \label{Mehta_seshadri} For a fixed parabolic type $\alpha_*$, there is a coarse moduli space $\mp:=\mp(r,\alpha_*,d)$ which is a projective irreducible normal variety,   parameterizing $\alpha_*$-semi-stable parabolic rank-r vector bundles of degree-d up to $S$-equivalence over the curve $C$.  Moreover,  the subspace $\mps \subset \mp$ of $\alpha_*$-stable parabolic bundles is an open smooth subset.  
\end{theorem}

For a line bundle $\delta \in \mathrm{Pic}^d(C)$,  we define $\smp_C(r,\alpha_*,\delta)=\{ E_* \in \mp_C(r,\alpha_*,d) \ / \ \det(E)\cong \delta\},$ the moduli space of $\alpha_*$-semistable parabolic bundles with determinant $\delta$,
 which is also projective irreducible normal variety.  
\subsection{Relative moduli spaces} \label{Relative moduli spaces}
In this subsection,  we will recall the existence of a relative version of the moduli spaces of semi-stable parabolic vector bundles over a family of smooth projective complex curves equipped with a family of parabolic divisors.    \\

Let $\pi_s: \mathcal{C} \longrightarrow S$ be a smooth family of projective curves of genus $g \ge 2$,  parameterized by an algebraic variety $S$ over $\mathbb{C}$ and let $
{\sigma_i : S \rightarrow \mathcal{C}}_{ i \in I}
$,  be $N$ section,   such that  $\forall i \neq j \in I$ and $\forall s \in S$,  we have:   $\sigma_i(s) \neq \sigma_j(s)$.  We denote by $D:= \sum_{i\in I} \sigma_i(S),$ the associated divisor (as the relative dimension of the map $\pi_s$ is one),  which will be seen as a family of parabolic degree $N$ divisors parameterized by the variety $S$ and let $\delta \in \mathrm{Pic}^d(\mathcal{C}/S)$\footnote{See $\cite{fantechi2005fundamental}$ for definition of relative Picard groups.}. \\

Let $\pi_e: \T \longrightarrow S$ be a $S$-variety.
 A relative family of $\alpha_*$-parabolic rank-r vector bundles of degree $d$ (resp.  determinant $\delta$) over $\mathcal{C}/S$ parameterized by $\T/S$ is a locally free sheaf $\eta$ over $\mathcal{C}\times_S \T$ together with the following data:
 \begin{itemize}
\item For  each $i \in I$,  we give a filtration of the vector bundle 
$ \mathcal{E}_{\sigma_i}:=\mathcal{E}\vert_{\sigma_i(S)\times_S \T}$
over $\sigma_i(S)\times_S \T\cong \T$ by subbundles  as follow
 $$ \mathcal{E}_{\sigma_i}=F^1_{i}(\eta) \supset  F^2_{i}(\eta) \supset \cdot\cdot\cdot \supset F^{\ell_i}_{i}(\eta) \supset F^{\ell_i+1}_{i}(\eta)=\{0\},$$ 
$$0 \leq a_1(i) < a_2(i)<...<a_{\ell_i}(i)<k,$$  
such that for each $j\in \{1,2,...,\ell_i\}$ we have:  
$\r \left( F^j_{i}(\mathcal{E})/F^{j+1}_{i}(\mathcal{E})\right)=m_j(i)$. Thus we get a parabolic structure over $\eta$, we denoted by $\eta_*$.
\item For each $t \in \T$ we set $\mathcal{C}_t:=\pi_s^{-1}\left( \pi_e(t)\right)$.  Then the vector bundle $\eta_* \vert_{\mathcal{C}_t}$ is a  $\alpha_*$-semistable parabolic bundle of degree $d$ respectively with determinant $\delta_t:=\delta \vert_{\mathcal{C}_t} \in \mathrm{Pic}^d(\mathcal{C}_t)$ with respect to the parabolic divisor $D_t:=\sum_{i\in I} \sigma_i(\pi_e(t)).$
 \end{itemize}
 Two relative families $\eta_*$ and $\eta'_*$ are said to be equivalent if there is a line bundle $L$ on $\T$ such that $\eta_* \cong \eta'_* \otimes \pi_e^*(L)$.   Hence, we define a functor 
 $$\begin{array}{ccccc}
  \underline{\mp_{\mathcal{C}}}:= \underline{\mp_{\mathcal{C}}(r,\alpha_*,d)} & : & S-schemes & \longrightarrow & Set \\
 & & \T& \longmapsto &   \underline{\mp_{\mathcal{C}}}(\T),    \\
\end{array}$$
which associates to a Noetherian $S$-scheme $\T$ the set of equivalent families of parabolic rank-r vector bundles over $\mathcal{C}/S$ parameterized by the scheme $\T/S$ of parabolic type $\alpha_*$ of degree $d$.  We define a subfunctor 
$$\underline{\smp_{\mathcal{C}}(r,\alpha_*,\delta)} \  \subset \ \underline{\mp_{\mathcal{C}}(r,\alpha_*,d)},$$ 
parameterizing  parabolic rank-r vector bundles over $\mathcal{C}/S$ of type $\alpha_*$ with determinant $\delta$.\\

Maruyama and Yokogawa constructed a relative version of the moduli space of semistable  parabolic vector bundles over a smooth family of projective curves in  \cite{yokogawa1993compactification,maruyama1992moduli} and $\cite{yokogawa1995infinitesimal}$.
\begin{theorem}[\cite{mehta1980moduli}] The functors defined above are representable by proper $S$-schemes denoted by
\begin{center}
$ \tilde{\pi_e}:  \mp_{\mathcal{C}/S}(r,\alpha_*,d)  \longrightarrow S, $ \quad and \quad $ \pi_e:  \smp_{\mathcal{C}/S}(r,\alpha_*,\delta)  \longrightarrow S.$
 \end{center} 
Their closed points parameterizes relative $S$-equivalence classes of rank-r semi-stable parabolic vector bundles of fixed type  $\alpha_*$ and degree $d$ respectively fix determinant $\delta$ over the family of marked curves $\mathcal{C}/S$.   
\end{theorem}

Let denote  $\chi^{par}:= \mathcal{C}\times_S \smp_{\mathcal{C}/S}(r,\alpha_*,\delta)$ the fibre product over $S$ (see diagram \eqref{fiber-product}). 
\begin{definition}[Universal family] A universal parabolic vector bundle over $\mathcal{C} \times_S \smp_{\mathcal{C}/S}(r,\alpha_*,d)$ is a family $\eta_*$ of parabolic vector bundle of rank-r with determinant $\delta$ of parabolic type $\alpha_*$ over the family of curves $\mathcal{C}/S$,  
such that: 
$  \forall \left[E_*\right] \in \smp_{\mathcal{C}/S}(r,\alpha_*,d),$
we have $\eta_* \vert_{\mathcal{C}_{E_*}}$ $S$-equivalent to $E_*$,
over the curve  $\mathcal{C}_{E_*}=\pi_s^{-1}\left( \pi_e \left( \left[E_*\right] \right) \right).$
\end{definition}
\begin{remark} 
\begin{enumerate}
\item A universal parabolic bundle if it exists is unique modulo equivalence of families.
\item In fact, existence of universal family is equivalent to  the isomorphism of functors
$$   \underline{\smp_{\mathcal{C}}(r,\alpha_*,d)}(-) \simeq \mathcal{H}om\left(-,\smp_{\mathcal{C}/S}(r,\alpha_*,d)\right),$$
in this case we say that the moduli space is a fine moduli space.
\end{enumerate}
\end{remark}
\begin{proposition}[$\cite{boden1999rationality} $, Proposition. 3.2] The moduli space of $\alpha_*$-parabolic-stable bundles is fine if and only if we have: \  $\gcd \{d,m_j(i) \vert i \in I,  1 \leq j \leq \ell_i \}=1$.
\end{proposition}
\begin{remark} As $\smp_{\mathcal{C}/S}(r,\alpha_*,\delta)$ is a a good quotient of a Hilbert quotient scheme, denoted $\mathcal{Z}^{ss}$, such that there is a universal bundle $\mathcal{E}_*$ on $\mathcal{C}\times_S\mathcal{Z}^{ss}$.  $\mathcal{E}_*$ may not descend to $\mathcal{C}\times_S \smp$, but objects such $\en^0(\eta)$,$\pe^0(\eta)$  and $\mathcal{A}^0(\eta)$ ect,  descend.   Recall that, a sheaf $\mathcal{F}$ on  $\mathcal{C}\times_S\mathcal{Z}^{ss}$ descends to $\mathcal{C}\times_S \smp$, if the action of scalar automorphisms of $\eta$ (relative to $\mathcal{Z}^{ss}$, on $\mathcal{F}$ is trivial.  This without confusion we pretend the existence of a universal bundle $\eta_*$ over $\mathcal{C}\times_S \smp_{\mathcal{C}/S}$ that we call virtual universal bundle.
\end{remark}
\section{Hecke modifications and filtered vector bundles} \label{Hecke Modification}
Let $E_* \longrightarrow C$ be a rank-r parabolic vector bundle of parabolic type $\alpha_*$ with respect to a parabolic divisor $D$ with determinant $\delta \in \mathrm{Pic}^d(C)$.  We associate  for all $i \in I$ and  $j \in \{1,2,3,...,\ell_i\}$ the following exact sequences 
$$0\longrightarrow \H_i^j(E) \hookrightarrow E \longrightarrow Q_i^j(E):= E_{x_i}/F_i^{j+1}(E)  \longrightarrow 0$$
the quotient sheaf $Q_i^j(E)$ is supported on $x_i$ of length $r_{j}(i)=\sum\limits_{q=1}^{j} m_q(i).$ The sub-sheaves $ \H_i^j(E) $ are locally free of rank-r and their determinants are given by:  
$ \delta_{j}(i):=\delta \otimes \mathcal{O}_C\left(-r_{j}(i)x_{i}\right)$,  
we denote their degree by  
$ d_j(i):=\deg \  \delta_j(i)=d-r_j(i),$
 and we set the integers
 $n_j(i)=\mathrm{gcd}(r,d_{j}(i))$ and $
 n=\mathrm{gcd}(r,d)$.
\begin{definition}[Hecke modifications]
We call the vector bundle $ \H_i^j(E)$ the Hecke modification of the parabolic bundle $E_*$ with respect to the subspace $F_i^{j+1}(E)\subset E_{x_i}$.   We set $\H_i^{0}(E)=E$.
\end{definition}
\begin{proposition}[Hecke filtrations]\label{Hecke-filtrations} Let $E_*$ be a parabolic rank-r vector bundle  with respect to the parabolic divisor $D$. Then, for each $i\in I$ the  Hecke modifications over $x_i \in D$,  satisfies for all $j \in  \{1,2,...,\ell_i\}$ the following inclusions
$$ E(-x_i)= \H^{\ell_i}_i(E)  \subset  \H^{\ell_i-1}_i(E)  \subset \cdot\cdot\cdot    \subset \H^2_i(E)   \subset  \H^{1}_i(E)  \subset \H^0_i(E)=E.$$
\end{proposition}
\bp We take the Hecke modifications over a point $x_i \in D$ for  $i\in I$ and let $j \in \{1,2,...,\ell_i \}$,  the $j$-th Hecke exact sequence
$$
\xymatrix{ 
0 \ar[r]  & \H^{j}_i(E) \ar[r] & E \ar[r]  &  Q_i^{j}(E) \ar[r]& 0, 
}
$$
where the last arrow is given by the composition 
$$
\xymatrix{ 
E \ar[rr]^{ev_{x_i}}  && E_{x_i}  \ar@{->>}[r] & Q_i^j(E)= E_{x_i}/F^{j+1}_i(E).
}
$$
 The inclusions 
$$ F^{j+1}_i(E) \supset F^{j+2}_i(E)$$
give a surjective maps 
$
Q_i^{j+1}(E) \twoheadrightarrow Q_i^j(E).
$
Then we get 
$$
\xymatrix{ 
0 \ar[r]  & \H^{j+1}_i(E) \ar[rr] \ar@{^{(}->}[rrd]^q && E  \ar[rr] \ar[d]^{id} & &  Q_i^{j+1}(E)    \ar[r] \ar@{->>}[d]&0 \\	
0 \ar[r] & \H^j_i(E) \ar[rr] && E  \ar[rr]_p&&  Q_i^j(E) \ar[r] & 0
}
$$
As the right diagram commutes  and the  map $p \circ q=0$,  we get that the image of the map $q$ is in the sub-sheaf $\H^j_i(E)$.  So as a conclusion, we get a filtration by rank-r locally free sub-sheaves 
$$ E(-x_i)= \H^{\ell_i}_i(E)  \subset  \H^{\ell_i-1}_i(E)  \subset \cdot\cdot\cdot    \subset \H^2_i(E)   \subset  \H^{1}_i(E)  \subset \H^0_i(E)=E.$$
\ep
\begin{remark} By the last proposition a rank-r  parabolic structure with respect to a parabolic divisor $D$ is equivalent to the following data: $(E,\H^*_*(E),\alpha_*)$ such that 
\begin{itemize}
\item $E$ a rank-r vector bundle over $C$. 
\item $\alpha_*=(k,\vec{a},\vec{m})$ is a parabolic type with respect to the divisor $D$.
\item for all $i \in I$,  we give a filtration by rank-r locally free subsheaves  
$$ E(-x_i)= \H^{\ell_i}_i(E)  \subset  \H^{\ell_i-1}_i(E)  \subset \cdot\cdot\cdot    \subset \H^2_i(E)   \subset  \H^{1}_i(E)  \subset \H^0_i(E)=E,$$
such that the torsion sheaves  $\H^{j}_i(E)/\H^{j+1}_i(E)$
 are supported at $x_i \in D$ and \\
$\mathrm{length}\left(\H^{j}_i(E)/\H^{j+1}_i(E)\right)=m_j(i).$
\end{itemize} 
\end{remark}
\paragraph{Classifying  maps}\label{section Classifying  maps}
Let $\mathcal{E}_*$ be a family of rank-r parabolic vector bundles of fixed parabolic type $\alpha_*$ over $(\mathcal{C},D)/S$ with fixed determinant $\delta \in \mathrm{Pic}^d(\mathcal{C}/S)$ parameterized by $\T/S$.  As the semi-stability is an open condition,  we get the following rational maps from $\T$:
\begin{itemize}
\item To the relative moduli space of parabolic semistable rank-r vector bundles of parabolic type $\alpha_*$
$$\begin{array}{ccccc}
 \psi_{\T} & : & \T & \dashrightarrow & \smp_{\mathcal{C}}(r,\alpha_*,\delta) \\ 
 & & t & \longmapsto & [ \mathcal{E}_{t_*} ]:=[ \eta_* \vert_{\mathcal{C}_t}] 
\end{array}$$
where $[\mathcal{E}_{t_*}]$ is the $S$-equivalence class of the semi-stable  parabolic bundle $\eta_{t_*}$.
\item  To the relative moduli spaces of the semistable rank-r vector bundles with fixed determinant $\delta_{j}(i) \in \mathrm{Pic}^{d_j(i)}(\mathcal{C}/S)$ for all $i \in  I$ and  $j \in \{1,2,...,\ell_i \}$ by associating  Hecke modifications (see subsection \ref{Hecke Modification})
$$\begin{array}{ccccc}
 \phi^{\T}_{i,j} & : & \T & \dashrightarrow & \mathcal{SU}_{\mathcal{C}}(r,\delta_{j}(i)) \\ 
 & & t & \longmapsto & \H^j_i(\mathcal{E}_t)
\end{array}$$
where $
\H^j_i(\mathcal{E}_t):=\mathrm{ker}\{\eta \longrightarrow Q^j_i(\eta)\}$ and $ 
\delta_j(i):=\delta \left(-r_j(i)\sigma_i(S)\right).$
\item The forgetful rational map (we forget the parabolic structure)
$$\begin{array}{ccccc}
\phi_{\T} &: & \T &  \dashrightarrow&  \mathcal{SU}_{\mathcal{C}}(r,\delta) \\ 
& &t & \longmapsto & \mathcal{E}_t 
\end{array}$$
\end{itemize}
We call these maps the classifying morphisms.
\subsection{Yokogawa-Maruyama point of  view}
In this subsection,  we give Yokogawa's point of view of  parabolic vector bundles and theirs moduli space.  Simpson $\cite{simpson1990harmonic}$ gives another description of parabolic vector bundles as filtered bundles,  which can be generalized to higher dimensions.   Maruyama and Yokogawa in \cite{maruyama1992moduli,yokogawa1991moduli,yokogawa1993compactification} give the construction of the relative moduli space of semistable parabolic vector bundles vector bundles using the new description and they prove that the two moduli spaces are isomorphic to each other as algebraic varieties.   \\

Let $C$ a smooth projective complex curve and $D= \sum_{i=1}^N x_i$ a reduced divisor over $C$.
\begin{definition} [Filtered vector bundles \cite{simpson1990harmonic} ]\label{Filtred-bundles}  A filtered rank-r bundle over the marked curve $(C,D)$ is a rank-r vector bundle $E$ over $C$ together with a filtrations  $E_{\bullet}=(E_{\lambda,i})_{i,\in I  \lambda \in \RR}$,  satisfying for all $i \in I$ the following conditions
\begin{enumerate}
\item Local freeness: $E_{\lambda,i}$ are locally free of rank-r,  $\forall \lambda \in \RR$  and $E_{0,i}=E$.
\item Decreasing: $E_{\lambda,i} \subset E_{\beta,i} $ \   for all  $ \lambda \geq \beta$.
\item Left continuous:  for $\varepsilon$ sufficiently small real number,  then $E_{\lambda-\varepsilon,i}=E_{\alpha,i}$.
\item Finiteness: the length of the filtration for $0 \leq \lambda \leq 1$ is finite.  
\item Periodicity: for all real number $\lambda$,  we have  $E_{\lambda+1,i}=E_{\lambda,i}(-x_i)$. 
\end{enumerate}
\end{definition}
\paragraph*{System of weights} Let $(E_{\lambda,i})_{i\in I,  \lambda \in \RR}$ be a filtered vector bundle with respect to the divisor $D$. Then we define the system of weights on $x_i$ for $i \in I$ as the ordered jumping real numbers in the reel interval $[0,1] $ i.e., $0 \leq \lambda \leq 1$ such that
 for $\varepsilon$ small enough we have:  $E_{\lambda,i} \neq E_{\lambda+\varepsilon,i}$. We will assume that  the jumping numbers are rational numbers. So we get for each $i \in I$ an ordered sequence of rational numbers: \ $   0 \leq \lambda_1(i) < \lambda_2(i) <...<\lambda_{\ell_i}(i)< 1,$
 where $\ell_i$ is the number of jumps at the point $x_i$. The $\mathrm{lenght} \left( \frac{E_{\lambda,i}}{E_{\lambda-1,i}} \right)$ is called the multiplicity of $\lambda_j(i)$.
\begin{remark}
We set $ E_{\lambda}=\bigcap\limits_{i=1}^N E_{\lambda,i}$,  we get a filtration $E_{\bullet}:=(E_{\lambda})_{\lambda \in \RR}$ that satisfies the first four points of the Definition \ref{Filtred-bundles}, and for the periodicity we get for each $\lambda \in \RR$, \ $E_{\lambda+1}=E_{\lambda}(-D).$ 
\end{remark}
\begin{definition}
Let $E_{\bullet}=(E_{\lambda})_{\lambda \in \RR}$ be a filtered rank-r bundle over the marked curve $(C,D)$. then we define 
\begin{enumerate}
\item Filtered degree:  
$\deg(E_{\bullet})=\int_0^1 \deg(E_{\lambda}) \mathrm{d}\lambda.$
\item Filtered slope:  $\mu(E_{\bullet})=\deg(E_{\bullet})/r.$
\end{enumerate}
\end{definition}
\begin{proposition}[Filtered bundles as Parabolic bundles]\label{Filtered bundles as Parabolic bundles}  Over a smooth marked curve $(C,D)$.
The notion of filtered rank-r vector bundles is equivalent to parabolic rank-r vector bundles with respect to the same divisor $D$ and stability conditions coincides.
\end{proposition}
\begin{remark}
$\cite{maruyama1992moduli}$ proved the following equality
$deg(E_{\bullet})=\pdeg(E_{\bullet})+\r(E) \ \deg(D).$
\end{remark}
 \paragraph{Classifying maps for filtered vector bundles:} Let $\eta_{\bullet}$ be a family of filtered rank-r bundles over the smooth family of marked curves $(\mathcal{C},D)$ over $S$ parameterized by a $S$-variety $\T$ of fixed determinant $\delta \in \mathrm{Pic}^d(\mathcal{C}/S)$ and fixed weights,  we get for each $\lambda \in \RR$ a rational map to the moduli space of semi-stable rank-r vector bundles of fixed determinant
 $$\begin{array}{ccccc}
\phi^{\T}_{\lambda} &: & \T &  \dashrightarrow&  \mathcal{SU}_{\mathcal{C}/S}(r,\delta(\lambda)) \\ 
& &t & \longmapsto & \eta_{\lambda}\vert_{\mathcal{C}_t} 
\end{array}$$
where for each $t \in \T$ we associate the curve $\mathcal{C}_t:=\pi_n^{-1}(t)=\pi_s^{-1}(\pi_e(t))$
and for each $\lambda \in \RR$ we associate the line bundle 
$\delta(\lambda):= \det(\eta_{\lambda}) \in \mathrm{Pic}^{d(\lambda)}(\mathcal{C}/S)$ of degree $d(\lambda)$ and denote $n(\lambda)=\gcd(r,d(\lambda))$.
\section{Line bundles over the moduli spaces of parabolic bundles}
In this section,  we recall the description of the Picard group of the relative moduli space of semi-stable vector bundles of fixed rank and determinant  and also its ample generator and its canonical line bundle.   \\

Let $\pi_s: \mathcal{C}\longrightarrow S$ be a smooth family of projective complex curves of genus $g \geq 2$ ($g \geq 3$ if the rank is 2 and the degree is even).  If the parabolic divisor is empty and that the parabolic type is trivial.  Note that the trivial parabolic structure is just the structure of a vector bundles and in this case parabolic semi-stability (resp. stability) coincide with  semistability (resp.  stability) of vector bundles.  Thus the relative moduli space $\smp_{\mathcal{C}}(r,0_*,d)$ in Theorem \ref{Mehta_seshadri} of rank-r parabolic bundles with determinant $\delta$ coincides with the coarse relative moduli space of semistable rank-r vector bundles with determinant $\delta$, denoted by 
$$ \mathcal{SU}_{\mathcal{C}/S}(r,d):=\smp_{\mathcal{C}/S}(r,0_*,d),$$
 is an irreducible normal variety over $S$.  
\subsection{Determinant line bundle}
Let $\eta$ be a family of semistable rank-r vector bundles with fixed determinant $\delta \in \mathrm{Pic}^d(\mathcal{C}/S)$ parametrized by a $S$-variety $\T$,  we get a cartesian diagram
$$
\def\cartesien{%
    \ar@{-}[]+R+<6pt,-1pt>;[]+RD+<6pt,-6pt>%
    \ar@{-}[]+D+<1pt,-6pt>;[]+RD+<6pt,-6pt>%
  }
\xymatrix{ 
 \mathcal{C}\times_S \T \ar[rr]^{p_n} \ar[d]_{ p_w} \cartesien & &   \T  \ar[d]^{p_e} &  \\ 
 \mathcal{C} \ar[rr]_{\pi_s=p_s} & &  S
}
$$
\begin{definition}[Determinant line bundle$\cite{knudsen1976projectivity}$]
 Let $\eta \rightarrow \mathcal{C}\times_S \T$ be a family of vector bundles.We define
$$\mathrm{det} R^{\bullet}p_{n_*}\left(\mathcal{E}\right):=\left(\det p_{n_*}(\mathcal{E})\right)^{-1} \otimes \det R^1p_{n_*}(\mathcal{E}),$$
which is an element of $\mathrm{Pic}(\T/S)$,  we call it the determinant line bundle associated to $\eta$ with respect to the map  $p_n:\mathcal{C}\times_S \T \longrightarrow \T$.
 \end{definition}

Drezet-Narasimhan gave the description of the ample generator of the relative Picard group $\mathrm{Pic}(\mathcal{SU}_{C}(r,\delta))$ and they describes the canonical bundle to the moduli space $\mathcal{SU}_{C}(r,\delta)$ for any line bundle $\delta$ over the curve $\mathcal{C}$. 
\begin{theorem} [ $\cite{drezet1989groupe}$, Theorems B \& F]\label{Drezet_Narasimhan} We have the following properties 
\item
\begin{enumerate}
\item The relative Picard group  $\mathrm{Pic}(\mathcal{SU}_{\mathcal{C}/S}(r,\delta)/S)$ is isomorphic to $ \ZZ \l$,  where $\l$ is an ample line bundle.
\item Set $n=\gcd(r,d)$, where $d=\deg(\delta)$.  Then the dualizing sheaf of $\mathcal{SU}_{\mathcal{C}/S}(r,\delta)$ is $K_{\mathcal{SU}_{\mathcal{C}/S}(r,\delta)/S} \cong \l^{-2n}.$
\end{enumerate}
\end{theorem}
\noindent Let $\eta$ be a  virtual universal bundle over $\mathcal{SU}_{\mathcal{C}/S}(r,\delta)$,  then
\begin{itemize}
\item  The relative ample generator of the Picard group is expressed as follow  $\l= \lambda(\eta \otimes p_w^*(F)).$ 
Where $F$ is a vector bundle over $\mathcal{C}$,  such that
$\r(F)=\frac{r}{n}$ and  $\deg(F)=-\frac{\chi(E)}{n}$.
\item The canonical bundle satisfy the equalities $\cite{laszlo1997line}$
$$\l^{-2n}=K_{\mathcal{SU}_{\mathcal{C}/S}(r,\delta)/S}=\lambda \left( \en^0(\eta) \right)^{-1}.$$
\end{itemize}

If $\T=\smp_{\mathcal{C}/S}(r,\alpha_*,\delta)$ and $k$ is large enough, the pullbacks under the classifying morphisms ${\phi_{i,j}},  \phi$ (we drop the reference to the parameter space) of the ample generators of $\mathrm{Pic} (\mathcal{SU}_{\mathcal{C}/S}(r,\delta_{j}(i))/S)$ and $\mathrm{Pic}( \mathcal{SU}_{\mathcal{C}/S}(r,\delta)/S)$, respectively, extend to all the space  $\smp_{\mathcal{C}/S}(r,\alpha_*,\delta)$.  We denote them by $\Theta_{j}(i)$ and $\Theta$, respectively.
\begin{theorem}[\cite{narasimhan1993factorisation}]\label{universal_pullback} Let $\eta$ be a relative family of rank-r vector bundle of fixed determinant $\delta \in \mathrm{Pic}^{d}(\mathcal{C}/S)$ parameterized by a $S$-scheme $\T$ over the family 
$ p_s : \mathcal{C} \longrightarrow S$,  then we have 
$$\phi^*_{\T} \left( \l \right)=\lambda(\eta)^{\frac{r}{n}}\otimes \det \left( \eta_{\sigma}\right)^{\aleph},$$
where 
$\phi_{\T}$ is the classifying morphism to $ \mathcal{SU}_{\mathcal{C}}(r,\delta)$ the moduli space of semi-stable rank-r bundles with determinant $\delta$,   $\sigma:S \longrightarrow \mathcal{C}$ any section of the map $p_s$,  and 
$$\aleph=\frac{d+r(1-g)}{n} \ \ \ \mathrm{and}  \ \ \ n=\gcd(r,d).$$
\end{theorem}
\subsection{Parabolic determinant line bundle}
Let $\mathcal{E}_*$ be a relative family of parabolic rank-r vector bundles of determinant $\delta \in \mathrm{Pic}^d(\mathcal{C}/S)$ and fixed parabolic type $\alpha_*$ over a smooth family of curves $(\mathcal{C},D)/S$ parameterized by $\T/S$.  Let $\pi_n: \mathcal{C} \times_S \T \longrightarrow \T$ the projection map. 

\quad Assume the following condition: \quad $kd+ \sum\limits_{i=1}^N \sum\limits_{j=1}^{\ell_i} m_j(i)a_j(i) \in r \ZZ \label{star}$.
\begin{definition} $\cite{BR93}$ \label{parabolic-determinant-by-Bisas}
We define the parabolic determinant line bundle as following 
$$ \lambda_{par}(\eta_*):= \lambda(\eta)^{k} \otimes \bigotimes\limits_{i=1}^N \bigotimes_{j=1}^{\ell_i} \left\lbrace \det \left( F^{j}_i(\mathcal{E})/F^{j+1}_i(\mathcal{E})\right)^{-a_j(i)}  \right\rbrace \otimes \det(\mathcal{E}_{\sigma})^{  \frac{k \chi_{par}}{r}},$$
which is a line bundle over $\T/S$,  where 
\begin{itemize}
\item $\mathcal{E}_{\sigma}:=\mathcal{E}\vert_{ \sigma(S) \times_S \T}$ for some section  $\sigma$ of the map  $\pi_s: \mathcal{C} \longrightarrow S$.
\item The determinant line bundle bundle:  $\lambda(\eta):=\det R^{\bullet}\pi_{n_*}(
 \eta):=\left( \det \pi_{n_*}\mathcal{E} \right)^{-1} \otimes \det R^1\pi_{n_*}(\eta).$
\item $\chi_{par}=d+r(1-g)+ \frac{1}{k}\sum\limits_{i=1}^N \sum\limits_{j=1}^{\ell_i} m_j(i)a_j(i) .$ 
\end{itemize}
\end{definition}
 \section{Hitchin connection in algebraic geometry} \label{chapter3}
 In this section, we introduce van Geemen-de Jong approach to the constructing connections over a push-forward of a line bundles by giving heat operators over the line bundles.  We will define connections,  heat operators,  the relation between them and give van Geemen-de Jong theorem,  which is an algebro-geometric analogues of Hitchin's theorem in K\"{a}hler geometry.   We follow \cite{van1998hitchin}. \\ 
 
In all the section,  we take $ \pi: \mathcal{M} \longrightarrow S,$ a smooth surjective morphism of regular  $\CC$-schemes,  we have the natural exact sequence associated to the differential map
\begin{equation}
\xymatrix{ 
0 \ar[r]  & T_{\mathcal{M}/S} \ar[r] &  T_{\mathcal{M}}  \ar[r]^{d\pi}  & \pi^* \left( T_S \right)  \ar[r] & 0. 
}
\label{dpi}
\end{equation}
Let $E$ be a locally free sheaf over $\mathcal{M}$.  We denote  $\mathcal{D}^{(q)}_{\mathcal{M}}(E)$  the sheaf of differential operators of order at most $q$ over $E$.  For each $ q \in \NN$,  we have a natural inclusion:    $ \mathcal{D}^{(q-1)}_{\mathcal{M}}(E) \hookrightarrow \mathcal{D}^{(q)}_{\mathcal{M}}(E).$ Hence we get the short exact sequence 
$$
\xymatrix{ 
0 \ar[r]  & \mathcal{D}^{(q-1)}_{\mathcal{M}}(E) \ar[r] &  \mathcal{D}^{(q)}_{\mathcal{M}}(E)  \ar[rr]^{\nabla_q} &&  \mathrm{Sym}^q(T_{\mathcal{M}})\otimes \en(E)  \ar[r] & 0,
}
$$
where $\mathrm{Sym}^q(T_{\mathcal{M}})$ is the $q$-th symmetric power and the natural  quotient map $\nabla_q$ is called the symbol map of order $q$. 
We define the sheaf $\mathcal{D}^{(q)}_{\mathcal{M}/S}(E)$ of relative differential operators with respect to the map $\pi: \mathcal{M} \rightarrow S $ as the sub-sheaf of operators that are $\pi^{-1}(\mathcal{O}_S)$-linear. 
By restriction to this sub-sheaf we get a map
$$ \nabla_q: \mathcal{D}^{(q)}_{\mathcal{M}/S}(E) \longrightarrow  \mathrm{Sym}^q(\T_{\mathcal{M}})\otimes \en(E)$$ 
with image in the sub-sheaf $\mathrm{Sym}^q( T_{\mathcal{M}/S}).$
Hence a relative symbol map
$$ \nabla_q: \mathcal{D}^{(q)}_{\mathcal{M}/S}(E) \longrightarrow  \mathrm{Sym}^q( T_{\mathcal{M}/S})\otimes \en(E) .$$ 
\subsection{Atiyah classes and connections on vector bundles}
We follow Atiyah's description of Atiyah algebroids and exact sequences $\cite{atiyah1957complex}$ in the context of vector bundles rather than working with principal bundles.
\begin{definition}[Atiyah Class] \label{Atiyah-definition} Let $E$ be a vector bundle over $\mathcal{M}$. Then the Atiyah exact sequence  associated to $E$ is given by the following  pull-back
$$
\xymatrix{ 
0 \ar[r]  & \en(E)  \ar[r] \ar@{=}[d] & \mathcal{A}_{\mathcal{M}}(E) \ar[r]^{\nabla_1} \ar@{^{(}->}[d] & T_{\mathcal{M}}  \ar[r] \ar@{^{(}->}[d]^{-\otimes id} & 0 \\	
0 \ar[r] & \en(E)  \ar[r]  & \mathcal{D}^{(1)}_{\mathcal{M}}(E)  \ar@{->>}[r]^{\nabla_1}  & T_{\mathcal{M}} \otimes \en(E)  \ar[r] & 0 }
$$
The sheaf $ \mathcal{A}_{\mathcal{M}}(E)$ is called the  Atiyah algebroid of $E$,  we denote its extension class by $at_{\mathcal{M}}(E)\in \mathrm{Ext}^1(T_{\mathcal{M}},\en(E)) \simeq \mathrm{H}^1(\mathcal{M},\Omega^1_{\mathcal{M}} \otimes \en(E))$\footnote{
 The isomorphism is because we deal with locally free sheaves. }. 
 For a line bundle $L$ over $\mathcal{M}$ the Atiyah sequence coincides with
\begin{equation}
\xymatrix{ 
  0 \ar[r] & \mathcal{O}_{\mathcal{M}} \ar[r] & \mathcal{D}^{(1)}_{\mathcal{M}}(L) \ar[r]^{\nabla_1}   & T_{\mathcal{M}}  \ar[r] & 0 }
  \label{atiyah-classe}
\end{equation}
\end{definition}

Note that the Atiyah class can be given by tensorize the Atiyah exact sequence \ref{Atiyah-definition} with the cotangent sheaf $\Omega^1_{\mathcal{M}}$
 $$
\xymatrix{ 
0 \ar[r]  & \en(E) \otimes \Omega^1_{\mathcal{M}} \ar[r] & \mathcal{A}_{\mathcal{M}}(E) \otimes \Omega^1_{\mathcal{M}} \ar[r]^{\nabla_1} & T_{\mathcal{M}}\otimes \Omega^1_{\mathcal{M}}  \ar[r]  & 0
}
$$
take the connecting morphism in the long exact sequence in cohomology 
$$\delta_1: \mathrm{H}^0 \left( \mathcal{M},\en( T_{\mathcal{M}}) \right) \longrightarrow \mathrm{H}^1 \left(\mathcal{M}, \en(E)\otimes \Omega^1_{\mathcal{M}} \right)$$
 then the class $at_{\mathcal{M}}(E)$ is given by $\delta_1(\mathrm{Id})$.  We have the following lemma in Atiyah \cite{atiyah1957complex}.
\begin{lemma} \label{multiplicity}  Let $X$ a smooth algebraic variety, $L$ a line bundle and $k$ a positive integer,  then we have an isomorphism of short exact sequences
  $$
\xymatrix{ 
0 \ar[r]  & \mathcal{O}_X \ar[r] \ar[d] & \mathcal{A}_X(L^k) \ar@{->>}[r]^{\nabla_1} \ar[d] & T_X  \ar[r] \ar[d]^{\mathrm{id}} & 0 \\	
0 \ar[r] &\mathcal{O}_X \ar[r]^{1/k}  & \mathcal{A}_X(L) \ar[r]^{\nabla_1}  & T_X  \ar[r] & 0 }
$$
 \end{lemma}

\quad  We define a relative version of Atiyah algebroid denoted by $\mathcal{A}_{\mathcal{M}/S}(E)$,  given by taking the pull-back 
$$
\xymatrix{ 
0 \ar[r]  & \en(E )\ar[r] \ar@{=}[d] & \mathcal{A}_{\mathcal{M}/S}(E)  \ar[r] \ar@{^{(}->}[d] & T_{\mathcal{M}/S}  \ar[r] \ar@{^{(}->}[d]^{\iota} & 0 \\	
0 \ar[r] & \en(E) \ar[r] & \mathcal{A}_{\mathcal{M}}(E)    \ar[r]^{\nabla_1}  & T_{\mathcal{M}}  \ar[r] & 0 }
$$
as an extension,  we have that  $ at_{\mathcal{M}/S}(E)$ in a section of $ \mathcal{E}xt^{1}(T_{\mathcal{M}/S},\en(E))\simeq R^1 \pi_*(\Omega^1_{\mathcal{M}/S} \otimes \en(E)).$ For a line bundle $L \in \mathrm{Pic}(\mathcal{M})$,  we denote this class by $[L] \in \mathrm{H}^0 \left( S, R^1\pi_*(\Omega^1_{\mathcal{M}/S}) \right)$. 
\\

For our purpose, we need the trace-free Atiyah algebroid of vector bundles with fix determinant.   We have a direct sum decomposition $\en(E)=\en^0(E) \oplus \mathcal{O}_{\mathcal{M}}$, and denote by $q: \en(E)\rightarrow \en^0(E)$ the first projection map.  Then, the trace-free Atiyah algebroid is given by the push-out of the standard Atiyah sequence by the map $q$ as follows
$$
\xymatrix{ 
0 \ar[r]  & \en(E)\ar[r] \ar[d]^q & \mathcal{A}_{\mathcal{M}}(E)  \ar[r]  \ar[d] & T_{\mathcal{M}}  \ar[r]  \ar[d] & 0 \\
0 \ar[r]  & \en^0(E)\ar[r] & \mathcal{A}^0_{\mathcal{M}}(E)  \ar[r]  & T_{\mathcal{M}}  \ar[r] & 0, 
 }
 $$
with the same method, we define the trace-free relative version Atiyah algebroid  $ \mathcal{A}_{\mathcal{M}/S}^0(E)$.
\subsection{Heat operators}
Let $L$ be a line bundle over $\mathcal{M}$ such that $\pi_*L$ is a locally free vector bundle over $S$.  We are interested in the subsheaf of the differential operators of degree 2 given by
$$ \mathcal{W}_{\mathcal{M}/S}(L):=   \mathcal{D}^{(1)}_{\mathcal{M}}(L)+\mathcal{D}^{(2)}_{\mathcal{M}/S}(L).$$
$\nabla_2$ is restriction of the symbol map to this sub-sheaf.
We define the subprincipal symbol 
$$ \sigma_S:  \mathcal{W}_{\mathcal{M}/S}(L) \longrightarrow  \pi^* T_S,$$
such that for $s \in L$ a local section of $L$ and $f$ a local section of $\mathcal{O}_S$ we have,  for all $D \in  \mathcal{W}_{\mathcal{M}/S}(L)$
$$\langle \sigma_S(H),d(\pi^*f) \rangle =d(\pi^*fs)-\pi^*fH(s).
$$
The elements of the sheaf $\mathcal{W}_{\mathcal{M}/S}(L) $ satisfies Leibniz rule (that follow from proprieties of the second symbol map) 
$$ H(fgs)= \left\langle \nabla_2(H),df \otimes dg \right\rangle s +fH(gs)+gH(fs)-fgH(s).
$$
Thus, we get a short exact sequence 
 \begin{equation}
 0 \longrightarrow  \mathcal{D}^{(1)}_{\mathcal{M}/S}(L) \longrightarrow
  \mathcal{W}_{\mathcal{M}/S}(L) \overset{\sigma_S \oplus \nabla_2}{\longrightarrow} \pi^* ( T_{S})\oplus \mathrm{Sym}^2(T_{\mathcal{M}/S} ) \longrightarrow 0.
  \label{heat-operator-sequence}
   \end{equation}
\begin{definition}[Heat operator $\cite{van1998hitchin}$] A heat operator $H$ on $L$  is a $\mathcal{O}_S$-linear map of coherent sheaves 
$$
H : T_S \longrightarrow \pi_*\mathcal{W}_{\mathcal{M}/S}(L)
$$
such that $\sigma_S \circ  \tilde{H}=Id,$ where $\tilde{H}$ is the $\mathcal{O}_{\mathcal{M}}$-linear map associate to $H$ by adjunction 
$$\tilde{H} : \pi^* T_S  \longrightarrow \mathcal{W}_{\mathcal{M}/S}(L).$$
\noindent A projective heat operator $H$ on $L$ is a heat operator with values in the sheaf $ \left( \pi_* \mathcal{W}_{\mathcal{M}/S}(L)\right)/\mathcal{O}_S $.
\end{definition}
\paragraph{Symbol of heat operators:} The symbol map of a (projective) heat operator $H$ is the map 
$$\rho_H:=\pi_*(\sigma_2) \circ H : T_S \longrightarrow \pi_* \mathrm{Sym}^2(T_{\mathcal{M}/S}).$$
\paragraph{The map $\mu_L$:} For any line bundle $L$ over  $\mathcal{M}$,  we associate the exact sequence  
$$ 0 \longrightarrow T_{\mathcal{M}/S} \longrightarrow \mathcal{D}_{\mathcal{M}/S}^{(2)} (L)/\mathcal{O}_{\mathcal{M}} \longrightarrow \mathrm{Sym}^2(T_{\mathcal{M}/S}) \longrightarrow 0,$$
and the first connector map in cohomology with respect to the map $\pi$ give rise to a map 
$$\mu_{L}: \pi_* \mathrm{Sym}^2( T_{\mathcal{M}/S}) \longrightarrow R^1 \pi_{*} \left( T_{\mathcal{M}/S} \right).$$

We have the following description.
\begin{proposition}[$\cite{welters1983polarized,pauly2023hitchin}$] \label{Welters} The first connecting morphism is given by the following formula: 
$ \mu_{L}=\cup \left[ L \right]-\cup \left( \frac{1}{2} \left[ K_{\mathcal{M}/S} \right] \right)$,  where $K_{\mathcal{M}/S}$ is the relative canonical line bundle.
\end{proposition}
\subsection{A heat operator for a candidate symbol}
As in Hitchin's theorem,  van Geemen-de Jong gave under what conditions  a candidate symbol  
$$\rho: T_S \longrightarrow \pi_* \left( \mathrm{Sym}^2(T_{\mathcal{M}/S})\right),$$
can be lifted to a (projective) heat operator, i.e. ,  is there a (projective) heat operator such that we have $\rho_H:=\sigma_S \circ H=\rho.$ The answer is given in the following theorem.
\begin{theorem}[van Geemen-de Jong, $\cite{van1998hitchin}$, \S 2.3.7]\label{van Geemen and De Jong}  Let $L \in \mathrm{Pic}(\mathcal{M})$ and $\pi: \mathcal{M} \longrightarrow S$ as before,  we have that if,   for a given map $\rho: T_S \longrightarrow \pi_* \mathrm{Sym}^2 T_{\mathcal{M}/S}$  
\begin{enumerate}
\item $\kappa_{\mathcal{M}/S}+\mu_{L}\circ \rho=0$,
\item The map  $ \cup [L]:\pi_* T_{\mathcal{M}/S} \longrightarrow R^1\pi_* \mathcal{O}_{\mathcal{M}}$,  is an isomorphism,  and
\item $\pi_* \mathcal{O}_{\mathcal{M}}=\mathcal{O}_S$.
\end{enumerate}
Then, there exists a unique projective heat operator $H$ whose symbol is $\rho$.
\end{theorem}
\begin{theorem}[Flatness criterion,  \cite{pauly2023hitchin} Theorem 3.5.1] \label{Flatness criterion} Under the  assumptions of the Theorem \ref{van Geemen and De Jong}, the projective connection  associated to the symbol $\rho$ is projectively flat if the following conditions holds.
\begin{enumerate}
\item  The symbol Poisson-commute with respect to the natural symplectic form over the relative cotangent bundle $T^{\vee}_{\mathcal{M}/S}$.i.e. \  For all local sections  $\theta, \theta'$ of $T_S$,  we have $ \{ \rho(\theta),\rho(\theta')\}_{T^{\vee}_{\mathcal{M}/S}}=0.$ 
\item The morphism $\mu_L$ is injective.
\item There are no vertical vector fields,  $\pi_* \left( T_{\mathcal{M}/S}\right)=0$.
\end{enumerate}
\end{theorem}
\section{Parabolic and strongly parabolic Atiyah sequences}
In this section, we prove the main theorem which  generalises $\cite{pauly2023hitchin}$ algebro-geometric construction of Hitchin connection over $\mathcal{SU}_{\mathcal{C}}(r,\mathcal{O}_{\mathcal{C}})$ the relative moduli space of rank-r vector bundles with trivial determinate over a smooth family of complex projective curves of genus $g \geq 2$,  to $\smp_{\mathcal{C}}(r,\alpha_*,\delta)$ the relative moduli space of parabolic rank-r vector bundles of fixed determinant $\delta \in \mathrm{Pic}(\mathcal{C}/S)$ and of fixed parabolic type $\alpha_*$.  \\

Let $S$ be a smooth complex algebraic variety,  we take a smooth family over  $S$ of projective marked curves $(\mathcal{C},D)$, where the divisor $D$ is given by $N$-section of the map $\pi_s$ such that the relative degree is $N$, let $\delta \in \mathrm{Pic}^d(\mathcal{C}/S)$ be a relative line bundles over the family of curves.  Let $\mathcal{U}$ be a family of rank-r vector bundles and fixed determinant $\delta$ and $\eta_*$ a family of  rank-r parabolic vector bundles of fixed parabolic type $\alpha_*$ and fixed determinant $\delta$ over $(\mathcal{C},D)/S$ parameterized by a $S$-schemes $\T$,  such that we have the following fibre product
$$
\xymatrix{ 
\mathcal{X}:=\mathcal{C} \times_S \T \ar[rr]^{\pi_n} \ar[d]_{\pi_w}   && \T  \ar[d]^{\pi_{e}} &  \\ 
\ \left( \mathcal{C},D \right) \ar[rr]_{\pi_s} && \ar@/^1pc/[ll]^{ \sigma_{i}} S
}
$$
Set $\mathcal{D}:=\pi_w^{-1}(D)=D \times_S \T$.  We define the quasi-parabolic and strongly quasi-parabolic Atiyah sequences and algebroids that we use to study deformation of marked curves equipped with quasi-parabolic vector bundles and to give an existence of Kodaira-Spencer map in the parabolic case . 
 \begin{definition}[Quasi-parabolic  Atiyah algebroid (QPA)]\label{QPA}
We take the pushout of the relative Atiyah exact sequence of the parabolic bundle $\eta_*$ by the inclusion $\en^0(\eta)\hookrightarrow \spe(\eta)^{ \vee}$,  we get 
$$
\xymatrix{ 
0 \ar[r]  & \en^0(\eta) \ar[r] \ar@{^{(}->}[d] & \mathcal{A}^{0}_{\mathcal{X}/ \T}(\eta) \ar[r] \ar[r] \ar[d] & T_{\mathcal{X}/\T} \ar[r]  \ar@{=}[d] & 0 \\	
0 \ar[r] & \spe^0(\eta )^{\vee} \ar[r]& \mathcal{A}_{1} \ar[r]   & \pi_w^* T_{\mathcal{C}/S} \ar[r] & 0}
$$
then the QPA sequence is given by tensorizing the exact sequence below by  $ \mathcal{O}_{\mathcal{X}}\left(\mathcal{-D} \right)$: $$
\xymatrix{ 
0 \ar[r] &  \pe^0(\eta)\ar[r] & \mathcal{A}^{0,par}_{\mathcal{X}/\T} (\eta)\ar[r]   & \pi_w^* T_{\mathcal{C}/S}(-D) \ar[r] & 0,
}
$$
and the QPA algebroid is given by: 
$ \mathcal{A}^{0,par}_{\mathcal{X}/\T} (\eta):=\mathcal{A}_{1}\otimes \mathcal{O}_{\mathcal{X}}\left(\mathcal{-D} \right).$
\end{definition}
\begin{definition} [Strongly quasi-parabolic Atiyah algebroid (SQPA)] \label{SQPA} We take the pushout of the  Atiyah exact sequence of the parabolic bundle $\eta_*$ by the inclusion $\en^0(\eta)\hookrightarrow \pe(\eta)^{ \vee} $,  we get 
$$
\xymatrix{ 
0 \ar[r]  & \en^0(\eta) \ar[r] \ar@{^{(}->}[d] & \mathcal{A}^{0}_{\mathcal{X}/\T}(\eta) \ar[r]  \ar[d] & T_{\mathcal{X}/\T} \ar[r] \ar@{=}[d]&0 \\	
0 \ar[r] &  \pe^0(\eta)^{\vee} \ar[r]& \mathcal{A}_{2} \ar[r]  & \pi_w^* T_{\mathcal{C}/S}  \ar[r] & 0
}
$$ 
then the SQPA sequence is given by  tensorizing exact sequence below by $ \mathcal{O}_{\mathcal{X}}\left(\mathcal{-D} \right)$:
 $$
\xymatrix{ 
0 \ar[r] & \spe^0(\eta)\ar[r] & \mathcal{A}^{0,par,St}_{\mathcal{X}/\T} (\eta)\ar[r]   &  \pi^*_w  T_{\mathcal{C}/S} \left(-D \right) \ar[r] & 0,
}
$$
and the SQPA algebroid is given by: $ \mathcal{A}^{0,par,St}_{\mathcal{X}/\T} (\eta):=\mathcal{A}_{2}\bigotimes \mathcal{O}_{\mathcal{X}} \left(\mathcal{-D} \right).$
\end{definition}
\section{Trace complexes theory over $\smp_{\mathcal{C}/S}(r,\alpha_*,\delta)$}
\subsection{Trace complexe theory over $\mathcal{SU}_{\mathcal{C}/S}(r,\delta)$} 
The main ingredient in \cite{pauly2023hitchin} 
is a description of the Atiyah class of the the determinant line bundle $\l$ over the moduli space $\mathcal{SU}_{\mathcal{C}/S}(r)$ by the Atiyah class of a universal family of rank-r vector bundles. In this subsection, we show the same relation on the moduli space $\mathcal{SU}_{\mathcal{C}/S}(r,\delta)$, based on  Beilinson-Schechtman's trace-complex theory, Sun-Tsai isomorphism (Theorem \ref{BEA} bellow) and Bloch-Esnault complex (Theorem \ref{BSBE}).  Here, we do not need the definition of the complex trace,  we use Sun-Tsai's characterization of the (-1)- Bloch-Esnault term as definition.   We recall the following fibre product 
$$
\xymatrix{ 
\mathcal{X}:=\mathcal{C}\times_S \mathcal{SU}_{\mathcal{C}/S}(r,\delta) \ar[rr]^{p_n} \ar[d]_{p_w}   & &  \mathcal{SU}_{\mathcal{C}/S}(r,\delta) \ar[d]^{p_e} &  \\ 
 \mathcal{C} \ar[rr]_{p_s= \pi_{s}} & &  S
}
$$

Let $\mathcal{U}$ be a universal vector bundle over  $\mathcal{C}\times_S \mathcal{SU}_{\mathcal{C}/S}(r,\delta)$.  The following theorem give a characterization of the (-1)-Bloch-Esnault algebra $^{0}\mathcal{B}^{-1}_{\mathcal{SU}_{\mathcal{C}/S}/S}(\mathcal{U}),$   that we will use as a definition. 
\begin{theorem} [\cite{sun2004hitchin}]\label{BEA} There is a canonical isomorphism  of short exact sequences 
$$
\xymatrix{ 
	0 \ar[r]  & T^{\vee}_{\mathcal{X}/\mathcal{SU}_{\mathcal{C}/S}} \ar[r] \ar[d]^{\cong} & \mathcal{A}^{0}_{\mathcal{X}/\mathcal{SU}_{\mathcal{C}}}(\mathcal{U})^{\vee} \ar[r] \ar[d]^{\cong} & \en(\mathcal{U})^{\vee}\ar[r] \ar[d]^{\cong}_{-Tr}&0 \\	
	0 \ar[r] & K_{\mathcal{X}/\mathcal{SU}_{\mathcal{C}/S}} \ar[r] & ^{0}\mathcal{B}^{-1}_{\mathcal{SU}_{\mathcal{C}/S}/S}(\mathcal{U}) \ar[r]     & \en(\mathcal{U})  \ar[r] & 0
}
$$
where $K_{\mathcal{X}/\mathcal{SU}_{\mathcal{C}/S}} $ is the relative canonical bundle with respect to the map $p_n$.
\end{theorem}
\begin{theorem} [Beilinson- Schechtman \cite{beilinson1988determinant}, Bloch-Esnault \cite{esnault2000determinant}]  \label{BSBE} There is a canonical  isomorphism of exact sequences over $ \mathcal{SU}_{\mathcal{C}/S}(r,\delta)$
 $$
\xymatrix{ 
	0 \ar[r]  &  R^1 p_{n_*} (K_{\mathcal{X}/\mathcal{SU}_{\mathcal{C}/S}}) \ar[r] \ar[d]_{2r \cdot \mathrm{id}}^{\cong} & R^1 p_{n_*} ( ^0\mathcal{B}^{-1}(\mathcal{U})) \ar[r] \ar[d]^{\cong} & {R^1 p_{n}}_* (\en^{0}(\mathcal{U})^{\vee})\ar[r] \ar[d]^{-Tr}_{\cong} &0 \\	
	0 \ar[r] & \mathcal{O}_{\mathcal{SU}_{\mathcal{C}/S}}\ar[r] & \mathcal{A}_{\mathcal{SU}_{\mathcal{C}/S}/S}\left(\lambda(\en^0 \left(\mathcal{U} \right)\right) \ar[r]     & T_{\mathcal{SU}_{\mathcal{C}/S}/S} \ar[r] & 0
}
$$
\end{theorem}

Combining this two results we get the following theorem,  proven in \cite{pauly2023hitchin} for $\delta=\mathcal{O}_{\mathcal{C}}$,  but their proof work for any relative line bundle $\delta \in \mathrm{Pic}^d(\mathcal{C}/S)$.
\begin{theorem}\label{BS-generalized}
There is a canonical  isomorphism of exact sequences over $ \mathcal{SU}_{\mathcal{C}/S}(r,\delta)$
$$
\xymatrix{ 
0\ar[r] & R^1 p_{n_*} (K_{\mathcal{X}/ \mathcal{SU}}) \ar[r] \ar[d]_{\frac{r}{n}\mathrm{Id}}^{\cong} & R^1 p_{n_*} \left( \mathcal{A}^0_{\mathcal{X}/\mathcal{SU}}(\mathcal{U})^{\vee} \right) \ar[r]  \ar[d]_{\cong} & R^1 p_{n_*} \left( \en^0(\mathcal{U})^{\vee} \right)  \ar[d]^{\mathrm{Id}}_{\cong}  \ar[r] & 0   \\	
0\ar[r] &  \mathcal{O}_{\mathcal{SU}_{\mathcal{C}}}\ar[r] & \mathcal{A}_{\mathcal{SU}_{\mathcal{C}}/S}(\mathcal{L}) \ar[r] ^{\nabla_1}     & T_{\mathcal{SU}_{\mathcal{C}}/S}  \ar[r] & 0 
}
$$ 
where $\l$ is the relative ample generator of the group $\mathrm{Pic}\left(\mathcal{SU}_{\mathcal{C}}(r,\delta)/S \right)$ and $n=\gcd(r,\deg(\delta))$.
\end{theorem}

\noindent \bp By Theorem \ref{BEA} and \ref{BSBE},  one has the following isomorphism of short exact sequences $\mathcal{SU}_{\mathcal{C}/S}(r,\delta)$
$$
\xymatrix{ 
0\ar[r] & R^1 p_{n_*}\left( K_{\mathcal{X}/\mathcal{SU}_{\mathcal{C}}} \right) \ar[r] \ar[d]_{2r\cdot \mathrm{id}}^{\cong} & R^1 p_{n_*} \left( \mathcal{A}^0_{\mathcal{X}/\mathcal{SU}_{\mathcal{C}}}(\mathcal{U})^{\vee} \right) \ar[r]  \ar[d]^{\cong} & R^1 p_{n_*} \left(\en^0(\mathcal{U})^{\vee} \right)  \ar[d]_{\cong}^{-Tr}   \ar[r] & 0    \\	
0\ar[r] &  \mathcal{O}_{\mathcal{SU}_{\mathcal{C}}}\ar[r]  & \mathcal{A}_{\mathcal{M}/S}\left(\lambda(\en^0(\mathcal{
U})) \right) \ar[r] ^{\nabla_1}     & T_{\mathcal{SU}_{\mathcal{C}}/S}  \ar[r] & 0 
}
$$
By Drezet-Narasimhan theorem \ref{Drezet_Narasimhan} and \cite{laszlo1997line},  we have  $\lambda(\en^0(\mathcal{
U}))=K_{\mathcal{SU}_{\mathcal{C}}}=\l^{-2n}.$
Hence, we  get the following isomorphism
$$
\xymatrix{ 
0\ar[r] & R^1 p_{n_*}\left( K_{\mathcal{X}/\mathcal{SU}_{\mathcal{C}}} \right) \ar[r] \ar[d]_{2r\cdot \mathrm{id}}^{\cong} & R^1 p_{n_*} \left( \mathcal{A}^0_{\mathcal{X}/\mathcal{SU}_{\mathcal{C}}}(\mathcal{U})^{\vee} \right) \ar[r]  \ar[d]^{\cong} & R^1 p_{n_*} \left(\en^0(\mathcal{U})^{\vee} \right)  \ar[d]_{\cong}^{-Tr}   \ar[r] & 0    \\	
0\ar[r] &  \mathcal{O}_{\mathcal{SU}_{\mathcal{C}}}\ar[r]  & \mathcal{A}_{\mathcal{M}/S}\left(\l^{-2n} \right) \ar[r] ^{\nabla_1}     & T_{\mathcal{SU}_{\mathcal{C}}/S}  \ar[r] & 0 
}
$$
 and by applying Lemma \ref{multiplicity} (for $k=-2n$ and $L=\l$),   we get 
 $$
\xymatrix{ 
0\ar[r] & R^1 p_{n_*}\left( K_{\mathcal{X}/\mathcal{SU}_{\mathcal{C}}} \right) \ar[r] \ar[d]_{2r\cdot \mathrm{id}}^{\cong} & R^1 p_{n_*} \left( \mathcal{A}^0_{\mathcal{X}/\mathcal{SU}_{\mathcal{C}}}(\mathcal{U})^{\vee} \right) \ar[r]  \ar[d]^{\cong} & R^1 p_{n_*} \left(\en^0(\mathcal{U})^{\vee} \right)  \ar[d]_{\cong}^{-Tr}   \ar[r] & 0    \\	
0\ar[r] &  \mathcal{O}_{\mathcal{SU}_{\mathcal{C}}}\ar[r]^{\frac{1}{2n}}  & \mathcal{A}_{\mathcal{M}/S}\left(\l \right) \ar[r] ^{-\nabla_1}     & T_{\mathcal{SU}_{\mathcal{C}}/S}  \ar[r] & 0 
}
$$
 The right vertical map is $-Tr$,  the vertical left map is $2r \mathrm{Id}$ and the extension class of the last exact map is $-2n \left[\l \right]$ in $\mathrm{H}^0\left(S, R^1 \pi_*\left(\Omega^1_{\mathcal{SU}_{\mathcal{C}/S}(r,\delta)/S}\right)\right)$.  Hence we conclude that the extension class of the exact sequence 
 $$
\xymatrix{ 
0 \ar[r] & R^1 p_{n_*}\left( K_{\mathcal{X}/\mathcal{SU}_{\mathcal{C}}} \right) \ar[r]  & R^1 p_{n_*} \left( \mathcal{A}^0_{\mathcal{X}/\mathcal{SU}_{\mathcal{C}}}(\mathcal{U})^{\vee} \right) \ar[r]   & R^1 p_{n_*} \left(\en^0(\mathcal{U})^{\vee} \right)  \ar[r] & 0   
}$$
 equals $\frac{n}{r} \left[\l\right]$.  This conclude the proof.
 \ep
 \subsection{Parabolic Bloch-Esnault complex }
Now,  we work over $\smp_{\mathcal{C}/S}:=\smp_{\mathcal{C}/S}(r,\alpha_*,\delta)$.  We denote by $\mathcal{E}_*$ a virtual universal parabolic bundle over $ \mathcal{X}^{par}=\mathcal{C}\times_S \smp_{\mathcal{C}/S}$.  For our need, we define the $(-1)$-term of the parabolic Bloch-Esnault complex as follow.
\begin{definition}
 We define the (-1)-term of the parabolic Bloch-Esnault $^0\mathcal{P}^{-1}(\eta)$,  as a pullback of the (-1)-term of the Bloch-Esnault complex $^0\mathcal{B}^{-1}(\eta)$, by the natural inclusion $\pe^0(\eta) \hookrightarrow \pe^0(\eta)$,  as follows 
$$
\xymatrix{ 
	0 \ar[r] & K_{\mathcal{X}^{par}/\smp} \ar[r] \ar@{=}[d] & ^0\mathcal{P}^{-1}(\eta) \ar[r] \ar[d]    & \pe^0(\eta) \ar[r]\ar[d]   & 0 \\
	0 \ar[r]  & K_{\mathcal{X}^{par}/\smp} \ar[r]  &  ^0\mathcal{B}^{-1}(\eta) \ar[r]  & \en^0(\eta) \ar[r] &0 
}
$$
where $ K_{\mathcal{X}^{par}/\smp}$ is the relative canonical line bundle relatively to the map $\pi_n$.
\end{definition}
 
\quad We apply $R^1\pi_{e_*}$  to the (-1)-Bloch-Esnault term exact sequence
$$
\xymatrix{ 
	0 \ar[r] &  \mathcal{O}_{\smp_{\mathcal{C}}} \ar[r] \ar@{=}[d] & R^1 \pi_{n_*} \left( ^0\mathcal{P}^{-1}(\eta) \right) \ar[r] \ar[d]    & R^1 \pi_{n_*} \left( \pe^{0}(\eta) \right) \simeq T_{\smp_{\mathcal{C}} / S}  \ar[r]\ar[d]   & 0 \\
	0 \ar[r]  &   \mathcal{O}_{\smp_{\mathcal{C}}} \ar[r]  &  R^1 \pi_{n_*}( ^0\mathcal{B}^{-1}(\eta)) \ar[r]  & R^1 \pi_{n_*}( \en^{0}(\mathcal{E})) \ar[r] &0 
}
$$
the exact sequence bellow is the pullback of the Bloch-Esnault exact sequence of the vector bundle $\eta$ seen as a family over the  space $\mathcal{SU}_{\mathcal{C}/S}(r,\delta)$ by the forgetful morphism map 
$$\phi: \smp_{\mathcal{C}/S}(r,\alpha_*,\delta) \longrightarrow \mathcal{SU}_{\mathcal{C}}(r,\delta)$$
which can be lifted to a map on the fibre product
$$\phi: \mathcal{C}\times_S\smp_{\mathcal{C}/S}(r,\alpha_*,\delta) \longrightarrow \mathcal{C}\times_S\mathcal{SU}_{\mathcal{C}/S}(r,\delta)$$

\noindent we can choose a virtual universal bundle $\mathcal{U}$ over  $ \mathcal{C}\times_S\mathcal{SU}_{\mathcal{C}/S}(r,\delta)$,  such that $\phi^*\left( \mathcal{U}\right) \cong \eta$.
Moreover, the differential map  $$d \phi: T_{\smp_{\mathcal{C}}/S} \longrightarrow \phi^*\left(T_{\mathcal{SU}_{\mathcal{C}}/S} \right)$$
is given by applying  $R^1 \pi_{n_*} $  to the natural inclusion $\pe^0 \left(\eta\right) \hookrightarrow \en^0\left( \phi^* \left(\mathcal{U}\right)\right).$\\

By construction we get an identification theorem in the parabolic configuration of Theorem \ref{BEA}.
\begin{proposition}\label{sun-tsai-parabolic} Let $\eta_*$ be a virtual universal parabolic bundle over $\smp_{\mathcal{C}/S}(r,\alpha_*,\delta)$.  Then, there is a canonical isomorphism $ ^0\mathcal{P}^{-1}(\mathcal{E})  \simeq \left[ \mathcal{A}^{0,par,St}_{\mathcal{X}^{par}/\smp} (\mathcal{E}) \left(\mathcal{D} \right)\right]^{\vee}$,  such that
$$
\xymatrix{ 
	0 \ar[r] &  K_{\mathcal{X}^{par}/\smp} \ar[r] \ar[d]^{\cong} & ^0\mathcal{P}^{-1}(\eta) \ar[r] \ar[d]^{\cong}   & \pe^0(\eta) \ar[r]\ar[d]^{\cong}   & 0 \\
	0 \ar[r]  &  K_{\mathcal{X}^{par}/\smp} \ar[r]  &  \left[ \mathcal{A}^{0,par,st}_{\mathcal{X}^{par}/\smp} (\eta) \left(\mathcal{D} \right) \right]^{\vee}   \ar[r]  & \pe^0(\eta) \ar[r] &0 
}
$$
\end{proposition}

Hence, we get the following parabolic version of \cite{pauly2023hitchin} Theorem 4.4.1.
\begin{theorem} \label{Parabolic BSBE}
Let $\eta_{\bullet}=(\eta_{\lambda})_{\lambda \in \RR}$ be a virtual universal filtered bundle over $\mathcal{M}_{\bullet}\simeq \smp_{\mathcal{C}/S}(r,\alpha_*,\delta)$. Then for each $\lambda \in \RR$, we have the following isomorphism of short exact sequences over $ \smp_{\mathcal{C}/S}(r,\alpha_*,\delta)$
$$
\xymatrix{ 
 R^1 \pi_{n_*} \left( K_{\mathcal{X}^{par}/\mathcal{M}{\bullet}} \right) \ar@{^{(}->}[r] \ar[d]^{\simeq}_{{\frac{r}{n(\lambda)}}} & R^1 \pi_{n_*} \left( \left[ \mathcal{A}^{0,par,st}_{\mathcal{X}^{par}/\M} (\eta_{\lambda} )(\mathcal{D} )\right]^{\vee} \right)\ar@{->>}[r]  \ar[d]^{\simeq} & R^1 \pi_{n_*} \left(\pe^0(\eta_{\lambda})\right) \ar[d]^{\cong}   \\	
 \mathcal{O}_{\M} \ar@{^{(}->}[r] & \mathcal{A}_{\M/S}\left( \Theta(\lambda) \right) \ar@{->>}[r] ^{\nabla_1}     & T_{\M/S} 
}
$$ 
where $\Theta(\lambda)$ is the pullback of the ample generator of the group $\mathrm{Pic}\left( \mathcal{SU}_{\mathcal{C}/S}(r,\delta_{\lambda})/S \right)$ by the classifying maps
 $$\begin{array}{ccccc}
\phi_{\lambda} &: & \mathcal{M}_{\bullet} & \longrightarrow &  \mathcal{SU}_{\mathcal{C}/S}(r,\delta(\lambda)) \\ 
& & \eta_{\bullet} & \longmapsto & \eta_{\lambda} 
\end{array}$$
set $d(\lambda)=\deg \ \delta(\lambda)$ and $n(\lambda)=\gcd \left(r, d(\lambda) \right)$,
which is equivalent to the equality $$\frac{r}{n(\lambda)} \Delta_{\lambda}= [ \Theta(\lambda)] \in \mathrm{H}^0\left( S,R^1 \pi_{e_*} \left( \Omega^1_{\M/S} \right)\right),$$
where we denote by $\Delta_{\lambda}$ the extension class of the first exact sequence.
\end{theorem}

This theorem is equivalent, in the parabolic setting, to the following theorem using Hecke modification.
 \begin{theorem}\label{Parabolic BSBE-2} Under the same hypothesis. Let $\eta_*$ be a virtual parabolic universal bundle,  we have the following isomorphism of short exact sequences 
 over $ \smp_{\mathcal{C}}(r,\alpha_*,\delta)$ 
$$
\xymatrix{ 
 R^1 \pi_{n_*} \left(K_{\mathcal{X}^{par}/\smp_{\mathcal{C}}} \right) \ar@{^{(}->}[r] \ar[d]^{\cong}_{{\frac{r}{n_j(i)}}} & R^1 \pi_{n_*} \left(\left[ \mathcal{A}^{0,par,st}_{\mathcal{X}^{par}/\smp_{\mathcal{C}}} (\H^j_i\left(\eta \right))(\mathcal{D})\right]^{\vee} \right)\ar@{->>}[r]  \ar[d]^{\cong} & R^1 \pi_{n_*} (\pe^0(\H^j_i\left(\eta \right))) \ar[d]^{\cong}  \\	
 \mathcal{O}_{\smp_{\mathcal{C}}} \ar@{^{(}->}[r] & \mathcal{A}_{\smp_{\mathcal{C}}/S}\left( \Theta_{j}(i) \right) \ar@{->>}[r] ^{\nabla_1}     & T_{\smp_{\mathcal{C}}/S} 
  }
$$ 
where $n_j(i):=\mathrm{gcd}(r,\deg (\delta_j(i)))$ as defined before.  We denote the extension class of the the first exact sequence by $\Delta_{j}(i)$ and the Atiyah class of a line bundle $L$ by $[L]$.  Then the theorem is equivalent to the equality of global sections   
$$ \frac{r}{n_j(i)} \Delta_{j}(i)= [ \Theta_{j}(i) ] \in \mathrm{H}^0\left(S,R^1 \pi_{e_*}(\Omega^1_{\mathcal{M}^{par}/S}) \right).$$
With the same hypothesis we have 
$$
\xymatrix{ 
 {R^1 \pi_{n}}_* \left(K_{\mathcal{X}^{par}/\smp_{\mathcal{C}}}\right) \ar@{^{(}->}[r] \ar[d]^{{\frac{r}{n}}}_{\cong }  & R^1 \pi_{n_*} \left( \left[\mathcal{A}^{0,par,st}_{\mathcal{X}^{par}/\smp_{\mathcal{C}}} (\mathcal{E})(\mathcal{D})\right]^{\vee}\right)  \ar@{->>}[r]  \ar[d]^{\cong} & {R^1 \pi_{n}}_* (\pe(\eta)) \ar[d]^{\cong}  \\	
 \mathcal{O}_{\smp_{\mathcal{C}}} \ar@{^{(}->}[r] & \mathcal{A}_{\smp_{\mathcal{C}}/S}(\Theta)  \ar@{->>}[r]^{\nabla_1}    & T_{\smp_{\mathcal{C}}/S}  
}
$$ 
which is equivalent to the equality $$\frac{r}{n} \Delta= [ \Theta] \in \mathrm{H}^0\left(S, R^1 \pi_{e_*}(\Omega^1_{\smp_{\mathcal{C}}/S})\right).$$
\end{theorem}
\bp  Modulo a shifting by a rational number $\lambda$ in the filtered configuration which correspondent to Hecke modifications in the parabolic sitting, it is sufficient to prove the theorem for $\lambda=0$. Take the forgetful map $ \phi : \smp_{\mathcal{C}/S}(r,\alpha_*,\delta) \rightarrow \mathcal{SU}_{\mathcal{C}/S}(r,\delta)$,
which can be lifted to the fibre product over $S$ and let $\eta_*$ be a universal parabolic bundle over $ \mathcal{C} \times_S \smp_{\mathcal{C}/S}$. Let denote by $\mathcal{U}$ a virtual universal bundle over $\mathcal{C}\times_S \mathcal{SU}_{\mathcal{C}/S}$ such that $\phi^* \left( \mathcal{U} \right) \cong \eta$. Take the pullback of the exact sequence given in Theorem \ref{BS-generalized} by the forgetful map $\phi$, we get 
$$
\xymatrix{ 
0\ar[r] & \phi^* \left( R^1 p_{n_*} (K_{\mathcal{X}/ \mathcal{SU}})\right) \ar[r] \ar[d]_{\frac{r}{n}\mathrm{Id}}^{\cong} & \phi^* \left( R^1 p_{n_*} \left( \mathcal{A}^0_{\mathcal{X}/\mathcal{SU}}(\mathcal{U})^{\vee} \right)\right) \ar[r] \ar[d]_{\cong} & \phi^* \left( R^1 p_{n_*} \left( \en^0(\mathcal{U})^{\vee} \right) \right) \ar[d]^{\mathrm{Id}}_{\cong} \ar[r] & 0 \\	
0\ar[r] & \mathcal{O}_{\smp}\ar[r] & \phi^*\left( \mathcal{A}_{\mathcal{SU}_{\mathcal{C}}/S}(\mathcal{L})\right) \ar[r] ^{\nabla_1} & \phi^* \left( T_{\mathcal{SU}_{\mathcal{C}}/S}\right) \ar[r] & 0 
}
$$ 
take the differential map
 $\mathrm{d} \phi: T_{\smp/S} \longrightarrow \phi^*\left( T_{\mathcal{SU}_{\mathcal{C}}/S} \right), $
 which correspond to taking the first direct image $R^1 \pi_{n_*}$ of the natural inclusion of sheaves 
 $$ \pe^0(\eta) \hookrightarrow \en^0(\eta)=\phi^* \left(\en^0(\mathcal{U}) \right)$$
  Now, we take the pull-back of this isomorphism of exact sequences by $\mathrm{d}\phi$ by the two realisations as follows
$$
\xymatrix{ 
0 \ar[r] & \mathcal{O}_{\smp_{\mathcal{C}}} \ar[r] \ar@{=}[d] & R^1 \pi_{n_*} ( ^0\mathcal{P}^{-1}(\eta)) \ar[r] \ar[d] & R^1 \pi_{n_*} (\pe^{0}(\eta)) \ar[r]\ar[d]^{d\phi} & 0 \\
0 \ar[r] & \mathcal{O}_{\smp_{\mathcal{C}}} \ar[r] \ar@{=}[d] & R^1 \pi_{n_*}( ^0\mathcal{B}^{-1}(\mathcal{E})) \ar[r] \ar@{=}[d] & R^1 \pi_{n_*}( \en^{0}(\eta)) \ar[r] \ar@{=}[d]&0 \\	
0\ar[r] & \phi^* \left( R^1 p_{n_*} (K_{\mathcal{X}/ \mathcal{SU}_{\mathcal{C}}})\right) \ar[r] \ar[d]_{\frac{r}{n}\mathrm{Id}}^{\cong} & \phi^* \left( R^1 p_{n_*} \left( \mathcal{A}^0_{\mathcal{X}/\mathcal{SU}_{\mathcal{C}}}(\mathcal{U})^{\vee} \right)\right) \ar[r] \ar[d]_{\cong} & \phi^* \left( R^1 p_{n_*} \left( \en^0(\mathcal{U})^{\vee} \right) \right) \ar[d]^{\mathrm{Id}}_{\cong} \ar[r] & 0 \\	
0\ar[r] & \mathcal{O}_{\smp_{\mathcal{C}}}\ar[r] \ar@{=}[d] & \phi^*\left( \mathcal{A}_{\mathcal{SU}_{\mathcal{C}}/S}(\mathcal{L})\right) \ar[r] ^{\nabla_1} & \phi^* \left( T_{\mathcal{SU}_{\mathcal{C}}/S}\right) \ar[r] & 0 \\
0 \ar[r] & \mathcal{O}_{\smp_{\mathcal{C}}} \ar[r] & \mathcal{A}_{\smp/S}(\phi^*(\l)) \ar[r] \ar[u] & T_{\smp/S} \ar[r] \ar[u]_{d\phi} & 0
}$$ 
By construction the first and the last exact sequences are isomorphic as they are pullbacks of isomorphic exact sequences by the differential map 
$$
\xymatrix{
 0 \ar[r] & R^1 \pi_{n_*} (K_{\mathcal{X}^{par}/\smp_{\mathcal{C}}}) \ar[r] \ar@{=}[d]_{\frac{r}{n} \mathrm{Id}} & R^1 \pi_{n_*} ( ^0\mathcal{P}^{-1}(\eta)) \ar[r] \ar[d]^{\cong} & R^1 \pi_{n_*} (\pe^{0}(\eta)) \ar[r] \ar[d]^{\cong} & 0 \\
 0 \ar[r] & \mathcal{O}_{\smp} \ar[r] & \mathcal{A}_{\smp/S}(\phi^*(\l)) \ar[r] & T_{\smp/S} \ar[r] & 0
 }
 $$
 where $\l$ is the ample generator of the Picard group of the space $\mathcal{SU}_{\mathcal{C}/S}(r,\delta)$ that we denote by $\Theta$. Note that we have the equalities 
 $$ K_{\mathcal{X}^{par}/\smp_{\mathcal{C}}} \cong \pi_w^*\left( K_{\mathcal{C}/S} \right) \cong \phi^* \left( K_{\mathcal{X}/ \mathcal{SU}_{\mathcal{C}}} \right) .$$
  We conclude the proof by applying Proposition \ref{sun-tsai-parabolic}, to obtain the following isomorphism of exact sequences
   $$
   \xymatrix{ 
   R^1 \pi_{n_*} (K_{\mathcal{X}^{par}/\smp_{\mathcal{C}}}) \ar@{^{(}->}[r] \ar[d]^{{\frac{r}{n}}}_{\cong } & R^1 \pi_{n_*} \left( \left[\mathcal{A}^{0,par,st}_{\mathcal{X}^{par}/\smp_{\mathcal{C}}} (\mathcal{E})(\mathcal{D})\right]^{\vee}\right) \ar@{->>}[r] \ar[d]^{\cong} & {R^1 \pi_{n}}_* (\pe(\eta)) \ar[d]^{\cong} \\	
   \mathcal{O}_{\smp_{\mathcal{C}}} \ar@{^{(}->}[r] & \mathcal{A}_{\smp_{\mathcal{C}}/S}(\Theta) \ar@{->>}[r]^{\nabla_1} & T_{\smp_{\mathcal{C}}/S}
    }
    $$ Hence, this conclude the proof. 
   \ep
\section{Parabolic Hitchin symbol map }
Let $\mathcal{E}_*$ be a virtual universal parabolic vector bundle of fixed parabolic type $\alpha_*$ over  $\mathcal{C} \times_S \smp_{\mathcal{C}/S}(r,\alpha_*,\delta)$.  We want to define a parabolic version of the Hitchin symbol map, as given in $ \cite{pauly2023hitchin}$ in section 4.3.  We will use the following notation:  $ \smp_{\mathcal{C}} := \smp_{\mathcal{C}/S}(r,\alpha_*,\delta)$. \\

\noindent \textbf{First approach:} We take the trace map which we denote by $B$ as follows
$$\small{
\xymatrix{ 
 \pi_{n_*} \left( \en^0(\eta) \otimes \pi_w^*  K_{\mathcal{C}/S} (D )  \right) \otimes \pi_{n_*} \left(\en^0(\eta) \otimes \pi_w^*  K_{\mathcal{C}/S} \left( D  \right) \right)  \ar[d]_B & (\phi,\psi) \ar@{|->}[d] \\ 
 \pi_{n_*}  \pi_w^* \left( K^{\otimes 2 }_{\mathcal{C}/S} \left( 2D \right) \right) & B(\phi,\psi)=\mathrm{Trace}(\phi \circ \psi)
}
}
$$
we take its restriction to the subbundle $ \spe^0(\eta) \subset  \en^0(\eta)$ so the left hand side is the cotangent bundle $T^{\vee}_{\smp_{\mathcal{C}}/S}$,  so we get 
$$
\xymatrix{ 
B : T^{\vee}_{\smp_{\mathcal{C}}/S} \otimes T^{\vee}_{\smp_{\mathcal{C}}/S}   \ar[r] & \pi_{n_*}  \pi_w^* \left( K^{\otimes 2 }_{\mathcal{C}/S} ( 2D ) \right)
}
$$
a simple calculation gives  $$\mathrm{Image} ( B ) \subset \pi_{n_*}  \pi_w^* \left( K^{\otimes 2}_{\mathcal{C}/S} ( D )  \right)$$
we denote by $B$ the following  restriction
$$
\xymatrix{ 
B : T^{\vee}_{\smp_{\mathcal{C}}/S} \otimes T^{\vee}_{\smp_{\mathcal{C}}/S}   \ar[r] & \pi_{n_*} \pi_w^* \left( K^{\otimes 2 }_{\mathcal{C}/S} (D) \right) 
}
$$
we dualize and by   Serre's  duality relative to  $\pi_{n}$ we get 
$$
B^{\vee} : \pi^*_{e} \left( R^1 \pi_{s_*} \left( T_{\mathcal{C}/S}\left( -D \right) \right) \right) \rightarrow  T_{\smp_{\mathcal{C}}/S}\otimes T_{\smp_{\mathcal{C}}/S} 
$$
\begin{definition} [Parabolic Hitchin Symbol]\label{def-Parabolic-Hitchin-Symbol-map} The parabolic Hitchin  symbol  $\rho^{par}$ is the  morphism given by 
$$
\rho_{par}:=\pi_{e_*}(B^{\vee}): R^1 \pi_{s_*} \left( T_{\mathcal{C}/S} ( -D) \right) \longrightarrow  \pi_{e_*}  \mathrm{Sym}^2 \left( T_{\smp_{\mathcal{C}}/S} \right) 
$$
\end{definition}
\noindent \textbf{Second approach:}  We consider the evaluation map of the  sheaf $\spe^{0}(\eta)\otimes \pi^*_w \left( \omega_{\mathcal{C}/S} ( D) \right)$  composed with the injection map  $\spe(\eta) \subset  \pe(\eta)$,  we get 
$$ \pi_n^* \ \pi_{n_*} \left( \spe^{0}(\eta)\otimes \pi^*_w  K_{\mathcal{C}/S} ( D)  \right) \overset{ev}{\longrightarrow} \pe^{0}(\eta)\otimes\pi^*_w \left( K_{\mathcal{C}/S}( D) \right)$$
we dualize
$$\pe^0(\eta)^{\vee} \otimes \pi^*_w \left( T_{\mathcal{C}/S} (-D) \right)\overset{ev^{\vee}}{\longrightarrow} \pi_n^* \left( \pi_{n_*} \left( \spe^0(\eta) \otimes \pi^*_w \ K_{\mathcal{C}/S} (D) \right)\right)^{\vee} $$
this morphism gives a map which  we denote by $ev^{\vee}$
$$\pi^*_w \left( T_{\mathcal{C}/S} (-D) \right) \overset{ev^{\vee}}{\longrightarrow} \pe^0(\eta) \otimes \pi_n^* \left( \pi_{n_*} \left( \spe^{0}(\eta) \otimes \pi^*_w \left( K_{\mathcal{C}/S} (D)  \right)  \right) \right)^{\vee}$$
by Serre's duality  relatively to $\pi_n$ 
$$\pi^*_w \left( T_{\mathcal{C}/S}(-D) \right) \overset{ev^{\vee}}{\longrightarrow} \pe^0(\eta)\otimes \pi_n^* \left( R^1 \pi_{n_*} \left( \pe^0(\eta) \right) \right)$$
we apply $\pi_{e_*} \circ R^1 \pi_{n_*}$  and by the projection formula,  we get
$$ \pi_{e_*}\left( R^1 \pi_{n_*}\left( ev^{\vee} \right) \right): R^1 \pi_{s_*} \left( T_{\mathcal{C}/S}(-D) \right) \longrightarrow \pi_{e_*}  \mathrm{Sym}^2 \left( T_{\smp_{\mathcal{C}}/S} \right).$$
\begin{lemma} This application coincide with the parabolic Hitchin  symbol  $\rho_{par}$. i.e., \quad  
$ \rho_{par}=\pi_{e_*}\left( R^1 \pi_{n_*} \left( ev^{\vee} \right) \right)
$.
\end{lemma}
\begin{proposition}\label{symbol-invariance} The symbol map $\rho_{par}$ is invariant under Hecke modifications.
\end{proposition} 

The proposition is a consequence of the following: take $E_*$ a parabolic vector bundle over a curve $C$ of parabolic type $\alpha_*$ with respect to a divisor $D$.  
Let $g \in \pe(E)$ be a parabolic endomorphism.  Then
\begin{lemma}
The trace is invariant under Hecke modifications.i.e.
\begin{center}
 $tr(\H_i^j(g))=tr(g)$ for all $i\in I$ and $j \in \{1,2,...,\ell_i \}$.
 \end{center} 
 \end{lemma}
 \bp For $i\in I$ and $j \in \{1,2,...,\ell_i \}$ take the Hecke modification of $E$ with respect to the subspace $F^{j+1}_i(E)$ so we get a sub sheaf $$ f: \H_i^j(E) \hookrightarrow E$$
 which is an isomorphism over $C\setminus\{x_i\}$,  thus  
\begin{center}
 $tr(\H_i^j(g))=tr(g)$  \quad   over \quad $C\setminus\{x_i\}$.
\end{center}
 The vector bundle $\H_i^j(E)$ inherited a parabolic structure, and $\H_i^j(g)$ is a parabolic endomorphism with respect to this  parabolic structure.
$$ 
\xymatrix{ 
E \ar[rr]^{g} && E  \\ 
 \H_i^j(E) \ar[rr]_{ \H_i^j(g)}  \ar[u]^f &&  \H_i^j(E) \ar[u]_f & 
}
$$
Now, let us describe the map $\H_i^j(g)_{x_i}:  \H_i^j(E)_{x_i} \rightarrow   \H_i^j(E)_{x_i} $.
 We have the decomposition of the map $g$ with respect to the quotient exact sequence
$$ 
\xymatrix{ 
0 \ar[rr] &&  F_i^{j+1}(E) \ar[rr] \ar[d]^{g\vert_{ F_i^{j+1}(E)}} &&   E_{x_i}  \ar[rr] \ar[d]^{g_{x_i}} && Q_i^j(E):= E_{x_i}/F_i^{j+1}(E)  \ar[rr] \ar[d]^{\overline{g}} && 0   \\ 
0 \ar[rr] &&  F_i^{j+1}(E) \ar[rr] &&   E_{x_i}  \ar[rr] & &  Q_i^j(E):=E_{x_i}/F_i^{j+1}(E)  \ar[rr]  &&0 
}
$$
thus we have
$$g_{x_i}=
\begin{pmatrix}
 g\vert_{ F_i^{j+1}(E)} &  \ast  \\
 0 & \overline{g}
\end{pmatrix} \Longrightarrow  tr(g_{x_i})=tr(g\vert_{ F_i^{j+1}(E)})+tr(\overline{g}).
$$
the Hecke modification $\H_i^j(E)$ fit in the same diagram  
$$ 
\xymatrix{ 
0 \ar[rr] &&  Q_i^j(E)  \ar[rr] \ar[d]^{\overline{g}} &&  \H_i^j(E)  \ar[rr] \ar[d]^{\H_i^j(g)_{x_i}} && F_i^{j+1}(E)  \ar[rr] \ar[d]^{g\vert_{ F_i^{j+1}(E)}} && 0   \\  
0 \ar[rr] && Q_i^j(E) \ar[rr] &&   \H_i^j(E)_{x_i}  \ar[rr] & & F_i^{j+1}(E)  \ar[rr]  &&0 
}
$$
hence we get 
$$\H_i^j(g)_{x_i}=
\begin{pmatrix}
  \overline{g}  & 0 \\
  \ast  & g\vert_{ F_i^{j+1}(E)}
\end{pmatrix} \Longrightarrow  tr(\H_i^j(g)_{x_i})=tr(\overline{g})+tr(g\vert_{ F_i^{j+1}(E)})= tr(g_{x_i}),
$$
There for one has globally the equality: $tr(g)=tr(\H_i^j(g)) \in \mathcal{O}_{C}.$  This ends the proof.
 \ep
\begin{proposition} \label{rho_isomorphism} The parabolic Hitchin symbol map $\rho_{par}$ is an isomorphism.
\end{proposition}
\bp Take the relative cotangent bundle over $\smp_{\mathcal{C}}$, denote it by  $ q:  T^{\vee}_{\smp_{\mathcal{C}}/S} \longrightarrow \smp_{\mathcal{C}}$,  the projection map. One get the following isomorphism  
$$ (\pi_e \circ q)_* \left( \mathcal{O}_{T^{\vee}_{\smp_{\mathcal{C}}/S}} \right) \cong \bigoplus\limits_{q \geq 0} \pi_{e_*} \mathrm{Sym}^q \left( T_{\smp/S} \right)  $$
and take the $\mathbb{G}_m$-action over the moduli space of the parabolic  Higgs bundles $\H iggs^P(\alpha_*)$ that contain the cotangent space $T^{\vee}_{\smp/S}$ as a big open subspace.  Thus elements of $ \pi_{e_*} \mathrm{Sym}^q \left( T_{\smp/S} \right) $ can be seen as regular functions over $ \mathcal{O}_{T^{\vee}_{\smp/S}}$ of degree 2 with respect to the action of $\mathbb{G}_m$,  that can be extend by Hartog's theorem to all the space $\H iggs^P(\alpha_*)$ of strongly-parabolic Higgs bundles.  As the parabolic Hitchin system is equivariant under the $\mathbb{G}_m$-action,  they are obtained  from the quadratic part of the parabolic Hitchin base given by the space 
$$\pi_{s_*} \left( K_{\mathcal{C}/S}^{\otimes2}\left(D \right) \right) \cong R^1 \pi_{s_*} \left( T_{\mathcal{C}/S}\left( -D \right) \right)$$ 
\ep
\section{Kodaira-Spencer map}

\subsection{Infinitesimal deformations of $(C,D,F^*_*(E))$}
In this subsecion, we study the infinitesimal deformations of $\mathbb{E}:=(C,D,E_*)$ a smooth marked projective curve of genus $g \geqslant 2$ and $D$ a reduced divisor of degree $N$ equipped with a quasi-parabolic rank-r vector bundle of fixed quasi-parabolic type $\vec{m}$.   
 \begin{theorem} \label{infinitesimal deformations} The infinitesimal deformations of  $\mathbb{E}=(C,D,E_*)$ are parametrized by $\mathrm{H}^1 \left( C, \mathcal{A}^{par}_{C} (E) \right)$.
\end{theorem}
\bp   Let $\mathcal{U}=\{U_{\lambda}\}_{\lambda}$  be an affine cover of the curve $C$ such that any open affine set contains at most one point of the  divisor $D$,  we set $U_{\lambda,\mu}=\mathrm{Spec}(A_{\lambda,\mu})$. 
Let $\mathbb{E}_{\varepsilon}$ be an infinitesimal deformations of $\mathbb{E}$ given by $(C_{\varepsilon},D_{\varepsilon}, {E_*}_{\varepsilon})$,  where the deformation $(C_{\varepsilon},D_{\varepsilon})$ is given by the 1-cocycle $\{\vartheta_{\lambda,\mu}\}$  in $\mathrm{H}^1(C,T_C(-D))$,  and the vector bundle $E_{\varepsilon}$ is given by the 1-cocycle $\{ \xi_{\lambda,\mu} \}$ with values in $\widehat{ \mathcal{A}_C(E)}$ given by the pull-back 
\begin{equation}
\xymatrix{ 
0 \ar[r] &  \en(E)  \ar[r] \ar[d] &  \mathcal{A}_C(E)  \ar[r] & T_C  \ar[r] & 0   \\  
0 \ar[r] &  \en(E)  \ar[r]  &  \widehat{\mathcal{A}_C(E)}  \ar[r] \ar@{^{(}->}[u]  & T_C(-D)  \ar[r]  \ar@{^{(}->}[u] & 0  
}
\label{the-first-subsheaf}
\end{equation} 
By definition,  a parabolic bundle is given for all $i \in I$ by the Hecke filtration (see Proposition \ref{Hecke-filtrations})
\begin{equation}
 \H_i^{\ell_i}(E)  \subset  \H_i^{\ell_i-1}(E)  \subset \cdot\cdot\cdot    \subset \H_i^2(E)   \subset  \H_i^{1}(E)  \subset  \H_i^0(E)=E.
\label{Filtration-in-the lemma}
\end{equation}
Hence the parabolic vector bundle ${E_*}_{\varepsilon}$ is given also for all $i \in I$ by filtrations of locally free sheaves
$$  \H_i^{\ell_i}(E_{\varepsilon})  \subset  \H_i^{\ell_i-1}(E_{\varepsilon})  \subset \cdot\cdot\cdot    \subset \H_i^2(E_{\varepsilon})   \subset  \H_i^{1}(E_{\varepsilon})  \subset  \H_i^0(E_{\varepsilon})=E_{\varepsilon},$$
and the 1-cocycle $\{\tau_{\lambda,\mu} \}$ must preserve this filtrations,  locally the sheaf $\H_i^{j}(E)$ is identified with a $A_{\lambda,\mu}[\varepsilon]$-submodule denoted by $M^{i,j}_{\lambda,\mu}[\varepsilon] \subset M^0_{\lambda,\mu}[\varepsilon]=M_{\lambda,\mu}[\varepsilon]$,  
$$   \newcommand{\incl}[1][r]
  {\ar@<-0.3pc>@{^(-}[#1] \ar@<+0.2pc>@{-}[#1]}
\xymatrix{ 
M_{\lambda,\mu}[\varepsilon]  \ar[rr]^{\tau_{\lambda,\mu}}  && M_{\lambda,\mu}[\varepsilon]  \\ 
M^{i,j}_{\lambda,\mu}[\varepsilon] \ar[rr] \incl[u]  && M^{i,j}_{\lambda,\mu}[\varepsilon] \incl[u] & 
}
$$
The fact that the diagram commutes is equivalent to the fact that the 1-cocycle $\{\xi_{\lambda,\mu} \}$ preserve the filtration given by the $A_{\lambda,\mu}$-modules $\{M^{i,j}_{\lambda,\mu}\}$ associated to the filtration \eqref{Filtration-in-the lemma}.  Hence the 1-cocycle $\{\xi_{\lambda,\mu} \}$ has values in the sheaf $ \widehat{\mathcal{A}^{par}_C(E)}$ defined as the subsheaf of $ \widehat{\mathcal{A}_C(E)}$ given locally by differential operators preserving the subsheaves $\H_i^j(E)$.  Hence,  the infinitesimal deformations of $\mathbb{E}=(C,D,E_*)$ are given by the cohomology group $\mathrm{H}^1(C,\widehat{\mathcal{A}^{par}_C(E)} )$.  Note that
the sheaf $\widehat{\mathcal{A}^{par}_C(E)}$ can be included in an exact sequence
$$ 
\xymatrix{ 
0 \ar[r] &  \pe(E)  \ar[r] &  \widehat{\mathcal{A}^{par}_C(E)}   \ar[r]^{\nabla_1}  & T_C(-D)  
}$$
where the map $\nabla_1$ is the restriction of the natural map $:\widehat{\mathcal{A}_C(E)} \rightarrow T_C(-D)$ given in the exact sequence \eqref{the-first-subsheaf}.
\\
\smallskip
To conclude the proof we need to show the following isomorphism  $\widehat{\mathcal{A}^{par}_C(E)} \cong\mathcal{A}^{par}_C(E) $.  Note that by definition of $\mathcal{A}^{par}_C(E)$ as push-out we have 
$$\mathcal{A}^{par}_C(E) :=\{ (f,\partial)\  / \ \  f \in \pe(E) \, \ \partial \in \mathcal{A}_C(E)(-D) \ \mathrm{and} \ (f,0) \sim (0,f) \ \mathrm{if} \ f \in \en(E)(-D)\}.$$
Thus we can define an $\mathcal{O}_C$-linear map $\varrho$ as follows
$$\begin{array}{ccccc}
\varrho  & : &  \mathcal{A}^{par}_C(E) & \longrightarrow &   \widehat{\mathcal{A}^{par}_C(E)}  \\
 & & (f, \partial)  & \longmapsto & f+\partial. 
\end{array}$$ 
Clearly the map $\varrho$ induces identity map on $\pe(E)$.  Let us prove that $\varrho$ is an isomorphism:
\begin{enumerate}
\item Injectivity:  Let $(f,\partial) \in \mathcal{A}^{par}_C(E)$ such that $\varrho(f,\partial)=0 \Leftrightarrow f+\partial=0 \Leftrightarrow \partial=-f$,  hence $\partial, f \in \en(E)(-D)$ by definition of $\mathcal{A}^{par}_C(E)$,  we have $(f,\partial)=(f,-f)\sim(f-f,0)=0 \in\mathcal{A}^{par}_C(E)$.    
\item Surjectivity: Let $\partial \in  \widehat{\mathcal{A}^{par}_C(E)}$ we associate its symbol $\nabla_1(\partial)\in T_C(-D)$.  Take a lifting  $\widehat{\nabla_1(\partial)} \in \mathcal{A}^{par}_C(E)$ (modulo $\pe(E)$),  which can be written $\widehat{\nabla_1(\partial)}=(f,\widehat{\partial})$,  where $\widehat{\partial} \in \mathcal{A}_C(E)(-D)$ with $\nabla_1(\widehat{\partial})=\nabla_1(\partial)$ and   $f$ any element in $\pe(E)$ .  Note that $\partial, \widehat{\partial} \in \mathcal{A}_C(E)(-D) \Rightarrow \partial-\widehat{\partial} \in \pe(E)$.  For $f=\partial-\widehat{\partial}$,  one has $q(\partial-\widehat{\partial},\widehat{\partial})=\partial$.
\end{enumerate}
Hence, we get an isomorphism of exact sequences 
$$
\xymatrix{ 
0 \ar[r] &  \pe(E)  \ar[r] \ar@{=}[d]^{\mathrm{Id}} & \mathcal{A}^{par}_C(E) \ar[r] \ar[d]_{\varrho}^{\cong} & T_C(-D)  \ar[r] \ar@{=}[d]^{\mathrm{Id}} & 0   \\  
0 \ar[r] &  \pe(E)  \ar[r]  &  \widehat{\mathcal{A}^{par}_C(E)}    \ar[r]  & T_C(-D)  \ar[r]  & 0  
}
$$
This concludes the proof. 
\ep
\begin{remark} Note that $ \cite{biswas2022infinitesimal}$ studied the infinitesimal deformations of $\mathbb{E}=(C,D,E_*)$.  Where their  definition of parabolic Atiyah algebroid $At(E_*)$,  coincides with the definition of the sheaf  $ \widehat{\mathcal{A}^{par}_C(E)}$ hence isomorphic to the parabolic Atiyah algebroid $\mathcal{A}^{par}_C(E)$.
 \end{remark}
\subsection{Parabolic Kodaira-Spencer map}
Let $\pi_s : (\mathcal{C},D) \longrightarrow S$ be a smooth family of projective marked curves parametrized by an algebraic variety $S$ and let $\pi_e :\smp_{\mathcal{C}}(r,\alpha_*,\delta) \longrightarrow S$ the relative moduli spaces of parabolic rank-r vector bundles of fixed parabolic type $\alpha_*$ with determinant $\delta \in \mathrm{Pic}^d(\mathcal{C}/S)$.  Let $\mathcal{E}_*$ a virtual universal parabolic vector bundle over $\mathcal{C}\times_S \smp_{\mathcal{C}}(r,\alpha_*,\delta)$ 
$$ \def\cartesien{%
    \ar@{-}[]+R+<6pt,-1pt>;[]+RD+<6pt,-6pt>%
    \ar@{-}[]+D+<1pt,-6pt>;[]+RD+<6pt,-6pt>%
  }
\xymatrix{ 
\mathcal{X}^{par} \ar[rr]^{\pi_{n}} \ar[d]_{\pi_{w}} \cartesien && \smp_{\mathcal{C}}(r,\alpha_*,\delta) \ar[d]^{\pi_{e}} &  \\ 
\ \left( \mathcal{C},D \right) \ar[rr]_{\pi_{s}} && \ar@/^1pc/[ll]^{ \sigma_{i}} S
}
$$
we use the following notation $\smp_{\mathcal{C}}:=\smp_{\mathcal{C}}(r,\alpha_*,\delta)$.  We have two fundamental maps 
\begin{itemize}
\item The Kodaira-Spencer of the family of marked curves: $\kappa_{\mathcal{C}/S}: T_S \longrightarrow R^1 \pi_{s_*} \left( T_{\mathcal{C}/S} (-D) \right),$\\
given as the first connecting morphism on cohomology of the short exact sequence 
$$0 \longrightarrow T_{\mathcal{C}/S}(-D) \longrightarrow \T_{\mathcal{C}}(-D) \longrightarrow \pi_s^* T_S \longrightarrow S.$$
\item The  Kodaira-Spencer of the family of moduli spaces:
$ \kappa_{\smp_{\mathcal{C}}/S}: T_S \longrightarrow R^1 \pi_{e_*} \left( T_{\smp_{\mathcal{C}}/S}  \right)$, 
where   $$T_{\smp/S}   \cong R^1 \pi_{n_*} \left( \pe^0(\eta)\right)$$
given as the first connecting morphism on cohomology of the short exact sequence 
$$0 \longrightarrow T_{\smp_{\mathcal{C}}/S} \longrightarrow T_{\smp_{\mathcal{C}}} \longrightarrow \pi_e^* T_S \longrightarrow S.$$
\end{itemize}
Take the QPA sequence of the bundle $\mathcal{E}_*$ over $ \mathcal{X}^{par}$
$$
\xymatrix{ 
0 \ar[r] & \pe^{0}(\mathcal{E})\ar[r] & \mathcal{A}^{0,par}_{\mathcal{X}^{par}/\mathcal{M}^{par}} (\mathcal{E}) \ar[r]    &  \pi^{*}_{w}\left( T_{\mathcal{C}} \left(-D \right) \right) \ar[r] & 0}
$$
As  $\pi_{e_*} \left( \pi^*_w\left( T_{\mathcal{C}} \left(-D \right) \right) \right) =0$ et $R^2 \pi_{n_*} \left( \pe^{0}(\eta)\right)=0$ (the relative dimension of $\pi_n$ is 1),
we apply $ R^1 \pi_{n_*}$ we get an  exact sequence on $\mathcal{M}^{par}$:
$$ \small {
\xymatrix{ 
0 \ar[r] &  T_{\mathcal{M}^{par}/S} \ar[r] & R^1 \pi_{n_*} \left(\mathcal{A}^{0,par}_{\mathcal{X}^{par}/\mathcal{M}^{par}} (\mathcal{E}) \right)\ar[r]    & R^1 \pi_{n_*} \left( \pi^{*}_{w} \left( T_{\mathcal{C}} \left(-D \right) \right) \right) \ar[r] & 0
}}
$$
The first connector  with respect to $\pi_{e_*}$ denoted $\Phi^{par}$,  and called  the parabolic Kodaira-Spencer map. 
\begin{proposition}\label{Kodaira-spencer decomposition} Then $\Phi^{par}$ commute with the Kodaira-Spencers of the two families: $\Phi_{par}\circ \kappa_{\mathcal{C}/S}=\kappa_{\smp_{\mathcal{C}}/S}$
\end{proposition}

\quad  Let $\partial=\cup \Delta$ is the first connector of the long exact sequence for $ \pi_e$ of the sequence  
$$
\xymatrix{ 
0 \ar[r] & \mathcal{O}_{\smp_{\mathcal{C}}}  \ar[r]& R^1 \pi_{n_*}\left( \left[\mathcal{A}^{0,par,st}_{\mathcal{X}^{par}/\smp_{\mathcal{C}}} (\eta)\left(\mathcal{D} \right)\right]^{\vee} \right) \ar[r]   &  T_{\smp_{\mathcal{C}}/S} \ar[r] & 0}
$$
given by applying $ R^1\pi_{n_*} $ to the dual of  the SQPA sequence tonsorized  by $\mathcal{O}_{\mathcal{\chi}^{par}}(\mathcal{D})$
$$
\xymatrix{ 
0 \ar[r] & \Omega^1_{\mathcal{\chi}^{par}/\smp_{\mathcal{C}}}  \ar[r] & \left[\mathcal{A}^{0,par,st}_{\mathcal{X}^{par}/\smp_{\mathcal{C}}} (\eta)\left(\mathcal{D} \right)\right] ^{\vee} \ar[r]    &  \pe^0(\eta)  \ar[r] & 0}
$$
\begin{theorem} [Parabolic version of proposition 4.7.1  \cite{pauly2023hitchin}] \label{parabolic version of proposition 4.7.1}
 The following diagram  commute
$$
\xymatrix{ 
R^1 \pi_{s_*} \left( T_{\mathcal{C}/\mathcal{S}} (-D )   \right)   \ar[rr]^{-\Phi_{par}} \ar[rd]_ {\rho_{par}}   & & R^1 \pi_{e_*} \left( T_{\smp_{\mathcal{C}}/S}  \right)& \\	
 & \pi_{e_*} \mathrm{Sym}^2 \left( T_{\smp_{\mathcal{C}}/S} \right)   \ar[ru]_{\partial}
} 
$$
i.e., $ \Phi_{par}+\partial \circ \rho_{par}=0.$
\end{theorem}
 \bp  We need  the following lemma 
\begin{lemma} [\cite{pauly2023hitchin} lemma 4.5.1, page 23.] \label{lemma} Let $X$ a scheme, $V$ and $L$ be, respectively, a vector and a line bundle on $X$. Moreover, let  $F \in \mathrm{Ext}^1 \left( L,V \right)$ 
$$\xymatrix{ 
0 \ar[r]  & V  \ar[r]^{\mathsf{i}} &  F  \ar[r]^{\pi}  & L  \ar[r]  &  0 
}
$$
 by taking the dual and tensorizing with $V \otimes L$,  we get 
$$ 0 \longrightarrow V \longrightarrow F^* \otimes V \otimes L \longrightarrow V^* \otimes V \otimes L \longrightarrow 0$$
consider the injection  
$$\begin{array}{ccccc}
 \psi & : & L & \longrightarrow & V^* \otimes V \otimes L   \\
 & & t & \longmapsto & Id_V\otimes t \\
\end{array}$$
 then there exist a canonical injection $\phi: F \longrightarrow F^* \otimes V \otimes L$ such that the following  diagram commutes 
$$\xymatrix{ 
0 \ar[r]  & V  \ar[r] \ar@{=}[d] &  F \ar@{^{(}->}[d]^{\phi} \ar[r]^{-\pi}  & L   \ar@{^{(}->}[d]^{\psi} \ar[r]  &  0 \\
 0  \ar[r]  &  V  \ar[r] & F^* \otimes V \otimes L \ar[r] & V^* \otimes V \otimes L  \ar[r] & 0   \\
 }
$$
\end{lemma}
Now,  we prove the proposition.  Take the parabolic Atiyah sequence on $\mathcal{X}^{par}$ of the universal bundle  $\mathcal{E}$ relative to $\pi_n$.  We  note : $\mathcal{A}^{par}:= \mathcal{A}^{0,par}_{\mathcal{X}^{par}/\mathcal{M}^{par}} \left( \mathcal{E}  \right)$ and $\mathcal{A}^{str}:= \mathcal{A}^{0,par,str}_{\mathcal{X}^{par}/\mathcal{M}^{par}} \left( \mathcal{E}  \right)$,  and take the evaluation map composed with the inclusion $ \spe(\mathcal{E}) \hookrightarrow  \pe(\mathcal{E}) $ 
$$\pi_n^* \ \pi_{n_*} \left( \spe^{0}(\eta)\otimes \pi^*_w  K_{\mathcal{C}/S}\left( D \right)  \right) \overset{ev}{\longrightarrow} \pe^{0}(\mathcal{E})\otimes \pi^*_w  K_{\mathcal{C}/S} \left( D \right)$$
We dualize
$$\pe^{0}(\eta)^{\vee} \otimes \pi^*_w \ T_{\mathcal{C}/S} \left( -D \right)\overset{ev^*}{\longrightarrow} \pi_n^*\ \pi_{n_*} \left( \spe^{0}(\eta) \otimes \pi^*_w K_{\mathcal{C}/S}\left( D \right)\right)^{\vee} 
$$
We get the following morphism of exact sequences
$$
\xymatrix{ 
	0 \ar[d]  & 0 \ar[d] \\
	\pe^{0}(\eta) \ar@{=}[r]\ar[d]  & \pe^{0}(\eta) \ar[d]  \\ 
	 \pe^{0}(\eta)\otimes  \mathcal{A}^{{par}^{\vee}} \otimes \pi^*_w T_{\mathcal{C}/S} \left(-D \right)   \ar[d] \ar[r]^{q} & \pe^{0}(\eta)\otimes \pi_n^*  \pi_{n_*}  \left( \mathcal{A}^{str} \otimes \pi^*_w K_{\mathcal{C}/S}\left( D \right)\right)^{\vee} \ar[d] \\ 
	\pe^{0}(\eta) \otimes \pe^{0}(\eta)^{\vee} \otimes  \pi^*_w  T_{\mathcal{C}/S} \left(-D \right)\ar[r]^{ev^{\vee}} \ar[d] & \pe^{0}(\eta) \otimes \pi_n^* \ \pi_{n_*} \left( \spe^{0}(\eta) \otimes \pi^*_w K_{\mathcal{C}/S}\left( D \right)\right)^{\vee} \ar[d] \\
	0  & 0 
}
$$
The map $q$  is given by taking the dual of the evaluation map 
$$ev:  \pi_n^*  \pi_{n_*}  \left( \mathcal{A}^{par} \otimes \pi^*_w K_{\mathcal{C}/S}\left( D \right)\right) \longrightarrow  \mathcal{A}^{par} \otimes \pi^*_w K_{\mathcal{C}/S} \left(D \right)$$
and compose it with the natural inclusion $  \mathcal{A}^{str} \hookrightarrow   \mathcal{A}^{par}$.
We apply lemma \ref{lemma} to the left exact sequence ( for $V= \pe^0(\eta)$, $L= \pi^*_w \ T_{\mathcal{C}/S} (-D )$ and $F= \mathcal{A}^{par} $),  and we apply the Serre duality relative to $\pi_n$ for the right exact sequence  we get the morphism of exact sequences
$$
\xymatrix{ 
	0 \ar[d]  & & & 0 \ar[d] \\
	\pe^{0}(\eta) \ar@{=}[rrr]\ar[d]  & & & \pe^{0}(\eta) \ar[d]  \\ 
	 \mathcal{A}^{par}  \ar[d] \ar[rrr] & & & \pe^{0}(\mathcal{E})\otimes \pi_n^* \ R^1 \pi_{n_*}  \left( \mathcal{A}^{{str}}(\mathcal{D})^{\vee} \right) \ar[d] \\ 
	 \pi^*_w \ T_{\mathcal{C}/S} (-D ) \ar[rrr] \ar[d] & & &  \pe^{0}(\eta) \otimes \pi_n^* \ R^1 \pi_{n_*} \left(  \pe^{0}(\eta)  \right) \ar[d] \\
	 0  & &  &0 
}
$$
The left exact sequence is the parabolic Atiyah sequence  where we multiply the map $  \mathcal{A}^{par}     \longrightarrow  \pi^*_w( \ T_{\mathcal{C}/S} \left(-D  \right)) $  by $-1$, as shown in Lemma  \ref{lemma}. \\

We apply  $R^1 \pi_{n_*}$, we get 
$$
\xymatrix{ 
	0 \ar[d]  & 0 \ar[d] \\
	T_{\smp_{\mathcal{C}}/S}\ar@{=}[r]\ar[d]  & T_{\smp_{\mathcal{C}}/S}\ar[d]  \\ 
	 R^1 \pi_{n_*} \left(  \mathcal{A}^{par} \right) \ar[d] \ar[r] &T_{\smp_{\mathcal{C}}/S} \otimes  R^1 \pi_{n_*}  \left( \mathcal{A}^{{str}}( \mathcal{D})^{\vee} \right) \ar[d] \\	 
	R^1 \pi_{n_*} \left(  \pi^*_w \left( T_{\mathcal{C}/S} \left(-D \right) \right) \right)\ar[r] \ar[d] & T_{\smp_{\mathcal{C}}/S}\otimes  T_{\smp_{\mathcal{C}}/S}\ar[d] \\
	0  & 0 
}
$$
The right exact sequence is the $R^1\pi_{e_*}$ applied to the dual of strongly parabolic Atiyah sequence tensorized by $T_{\smp_{\mathcal{C}}/S}$,  and the first connecting homomorphism in cohomology with respect to $ \pi_{e}$ is given by cup product with the class $\Delta$. Take the first connecting homomorphism  of the long exact sequence with respect to the map $ \pi_e$,  we get 
$$
\xymatrix{   \ar @{} [ddrr] |{\square}
\pi_{e_*}  R^1 \pi_{n_*} \left(  \pi^*_w \left(T_{\mathcal{C}/S} \left(-D \right) \right)\right)  \ar[rr]  \ar[dd]  & &  \pi_{e_*} \left(  T_{\smp_{\mathcal{C}}/S} \otimes  T_{\smp_{\mathcal{C}}/S} \right) \ar[dd]^{\cup \Delta} \\	 \\
	R^1 \pi_{e_*} \left( T_{\smp_{\mathcal{C}}/S} \right) \ar@{=}[rr] &  & R^1 \pi_{e_*}  \left( T_{\smp_{\mathcal{C}}/S} \right)
	}
$$
we have   $$\pi_{e_*}  R^1 \pi_{n_*} \left(  \pi^*_w \left( T_{\mathcal{C}/S} \left(-D \right) \right)\right) \simeq R^1 \pi_{s_*} \left( T_{\mathcal{C}/S} \left(-D \right) \right) $$
 we get, the following commutative diagram
$$
\xymatrix{  \ar @{} [ddrr] |{\square}
R^1 \pi_{s_*} \left( \mathcal{T}_{\mathcal{C}/S} \left(-D \right) \right) \ar[rr]^ {\rho_{par}}   \ar[dd]_{ - \Phi^{par}} & & \pi_{e_*}  \left(  T_{\smp_{\mathcal{C}}/S} \otimes  T_{\smp_{\mathcal{C}}/S} \right) \ar[dd]^{ \partial} \\ \\	
	R^1 \pi_{e_*} \left( T_{\smp_{\mathcal{C}}/S} \right) \ar@{=}[rr] & & R^1 \pi_{e_*} \left( T_{\smp_{\mathcal{C}}/S} \right)
	}
$$
Thus conclude the proof.
\ep
\subsection{Some equalities and consequences}We recall the equalities given in Theorem \ref{Parabolic BSBE-2}. For all $i \in I$ and $j \in \{1,2,...,\ell_i\}$ one has
\begin{enumerate} 
\item $
 \frac{r}{n} \Delta= [ \Theta] \in \mathrm{H}^0\left(S, R^1 \pi_{e_*} (\Omega^1_{\smp/S})\right)$, and
\item $ \frac{r}{n_j(i)} \Delta_{j}(i)= [ \Theta_{j}(i) ] \in \mathrm{H}^0\left(S, R^1 \pi_{e_*} (\Omega^1_{\smp/S})\right)$.
\end{enumerate}
We denote the associated applications given by the contraction with the classes $\Delta$ and $\Delta_{j}(i)$, respectiverly by $$\partial := \cup \Delta : \pi_{e_*} \left( \mathrm{Sym}^2 \left( T_{\smp/S} \right) \right) \longrightarrow R^1 \pi_{e_*} \left( T_{\smp/\mathcal{S}} \right),$$
and $$\partial_j(i):= \cup \Delta_j(i): \pi_{e_*} \left( \mathrm{Sym}^2 \left( T_{\mathcal{M}^{par}/S} \right) \right) \longrightarrow R^1 \pi_{e_*} \left( T_{\mathcal{M}^{par}/\mathcal{S}} \right).$$

Combining the above equalities, we get the following result.
\begin{theorem} \label{equalities} Assume that the family $(\pi_s: (\mathcal{C},D) \rightarrow S)$ is versal \footnote{The Kodaira-Spencer of the family of marked curves is an isomorphism}. Then for all $i \in I$ and $j \in \{1,2,...,\ell_i \}$ we have the equalities over the moduli space $\smp_{\mathcal{C}/S}(r,\alpha_*,\delta)$
\begin{center}
 $\cup [\Theta]\circ \rho_{par}=- \frac{r}{n}\cdot \Phi_{par} $ \quad and  \quad
 $\cup [\Theta_j(i)]\circ \rho_{par}=- \frac{r}{n_j(i)} \cdot \Phi_{par} $.
\end{center}
\end{theorem}
\bp The first equality is a direct consequence of Proposition \ref{parabolic version of proposition 4.7.1}, where we have 
$ \partial \circ \rho_{par}= -\Phi_{par},$
we multiply the equality by $\frac{r}{n}$ and use the first equality above
$$\frac{r}{n} \partial \circ \rho_{par}=\cup [\Theta]\circ \rho_{par}= -\frac{r}{n}\Phi_{par}.$$
Let fix $i \in I$ and $j \in \{1,2,...,\ell_i \}$. We take the Hecke isomorphism over $S$
$$\def\commutatif{\ar@{}[u]|{\circlearrowleft}}\xymatrix{ \H^j_i : \smp_{\mathcal{C}/S}(r,\alpha_*,\delta) \ar[rr]^{\cong} \ar[dr]_{\pi_e} & & \smp_{\mathcal{C}/S} \left (r,\H_i^j(\alpha_*),\H_i^j(\delta)\right) \ar[dl]^{\pi_e^{i,j}} \\&	S \commutatif 
	}$$
where the map $\H^j_i $ is given by Hecke modification equipped with its natural parabolic structure.  Then by Proposition \ref{Kodaira-spencer decomposition} applied over $\smp:=\smp_{\mathcal{C}/S}(r,\alpha_*,\delta)$ and $\smp_{i,j}:=\smp_{\mathcal{C}/S} \left (r,\H_i^j(\alpha_*),\H_i^j(\delta) \right) $, the following diagram commute under the assumption that the map $\kappa_{\mathcal{C}/S}$ is an isomorphism
$$
\xymatrix{ 
&&&& R^1\pi_{e_*}\left( T_{\smp/S}\right) \ar[dd]_{\cong}^{\H_i^j} \\
T_S \ar[rrr]_{\cong}^{\kappa_{\mathcal{C}/S}} \ar@/^1pc/[rrrru]^{\kappa_{\pi_e}} \ar@/_1pc/[rrrrd]_{\kappa_{\pi^{i,j}_e}} &&& R^1\pi_{s_*} \left( T_{\mathcal{C}/S}(-D) \right) \ar[ru]^{\Phi_{par}} \ar[rd]_{\Phi^{i,j}_{par}} \\
&&&& R^1\pi^{i,j}_{e_*}\left( T_{\smp_{i,j}/S}\right)
	}$$
	In fact by Proposition \ref{Kodaira-spencer decomposition} applied over $\smp$ and $\smp_{i,j}$, one has 
\begin{equation}
 \H_i^j \circ \Phi_{par}= \Phi^{i,j}_{par}. 
\label{Kodaira-Spencer-Hecke}
\end{equation}
Now, we define parabolic Hitchin symbol map $\rho^{i,j}_{par}$ over the moduli space $\smp_{i,j}$ (see definition \ref{def-Parabolic-Hitchin-Symbol-map}). We have the following commutative diagram
$$\def\commutatif{\ar@{}[rrdd]|{\circlearrowleft}}
 \xymatrix{
\commutatif && \pi_{e_*} \left( \mathrm{Sym}^2 T_{\smp/S} \right) \ar[rr]^{\cup [\Theta_j(i)]} \ar[dd]_{\cong}^{\H_i^j} \commutatif && R^1\pi_{e_*}\left( T_{\smp/S}\right) \ar[dd]_{\cong}^{\H_i^j} \\
   R^1\pi_{s_*} \left( T_{\mathcal{C}/S}(-D) \right) \ar@/^1pc/[rru]^{\rho_{par}} \ar@/_1pc/[rrd]_{\rho^{i,j}_{par}} \
   \\
 && \pi_{e_*} \left( \mathrm{Sym}^2 T_{\smp_{i,j}/S} \right) \ar[rr]_{\cup [\Theta_j(i)]} && R^1\pi^{i,j}_{e_*}\left( T_{\smp_{i,j}/S} \right) 
  }$$
The first diagram commute by Proposition \ref{symbol-invariance}. Hence, by the above diagram one has
\begin{align*}
\cup [\Theta_j(i)] \circ \rho_{par} &=\left( (\H_i^j)^{-1} \circ \cup [\Theta_j(i)] \circ \H_i^j \right) \circ \rho_{par}= (\H_i^j)^{-1} \circ \left( \cup [\Theta_j(i)] \circ \rho^{i,j}_{par} \right) 
\end{align*}
we apply Proposition \ref{parabolic version of proposition 4.7.1}, we get 
\begin{align*}
\cup [\Theta_j(i)] \circ \rho_{par} &=(\H_i^j)^{-1} \circ \left( -\frac{r}{n_j(i)} \Phi_{par}^{i,j} \right)=-\frac{r}{n_j(i)} \Phi_{par},
\end{align*}
the last equality is given by equation \eqref{Kodaira-Spencer-Hecke}.
This conclude the proof.
\ep
\section{Some line bundles over $\smp_{\mathcal{C}/S}(r,\alpha_*,\delta)$}
\paragraph{Parabolic determinant bundle:}  Let $\mathcal{E}_*$ be a family of parabolic rank r vector bundle of fixed parabolic type $\alpha_*$ over a smooth  family of curves $\mathcal{C}/S$ parametrized by a $S$-variety $\T$.   Let $p : \mathcal{C} \times_S \T \longrightarrow \T$ the projection map.  We recall the definition of the parabolic determinant line bundle under the hypothesis \eqref{star}
$$ \lambda_{par}(\eta_*):= \lambda(\mathcal{E})^{k} \otimes \left( \bigotimes\limits_{i=1}^N \bigotimes_{j=1}^{\ell_i} \left\lbrace \det \left( F^{j}_i(\eta)/F^{j+1}_i(\eta)\right) ^{ -a_j(i)} \right\rbrace \right) \otimes \det(\eta_{\sigma})^{\frac{k \chi_{par}}{r}}$$
which is a line bundle over $\T$,  see definition \ref{parabolic-determinant-by-Bisas}.  Pauly in \cite{pauly1996espaces} give another definition as following 
\begin{definition}[Parabolic determinant bundle]\label{parabolic_determinant_2} Let $\mathcal{E}_*$ be a family of parabolic rank-r vector bundles of parabolic type $\alpha_*$ over a smooth family of curves $\pi_s: \mathcal{C} \longrightarrow S$ parameterized by a $S$-variety $\T$,  then we have  
$$ \Theta_{par}(\eta_*):= \lambda(\eta)^{k} \otimes \left( \bigotimes\limits_{i=1}^N \bigotimes_{j=1}^{\ell_i-1} \left\lbrace \det \left( \eta_{\sigma_i}/F^{j+1}_i(\eta)\right)^{ p_{j}(i)} \right\rbrace \right) \otimes \det(\mathcal{E}_{\sigma})^{e} $$
where the determinant is with respect to the projection $\mathcal{C} \times_S \T \longrightarrow \T$ and for all $i \in I$ and $j \in \{1,2,...,\ell_i-1\}$
\begin{itemize} 
\item $\eta_{\sigma_i}:=\eta\vert_{\sigma_i(S)\times_S \T}$,  where $\sigma_i:S \longrightarrow \mathcal{C}$ the parabolic section of $\pi_s$.
\item $p_j(i)=a_{j+1}(i)-a_j(i)$  and $r_j(i):=\sum\limits_{i=1}^{q} m_i(q)=\dim(\eta_{\sigma_i}/F^{j+1}_i(\eta))$.
\item $re= k\chi -\sum\limits_{i=1}^N \sum\limits _{j=1}^{\ell_i-1} p_j(i) r_j(i)$,  where $ \chi=d+ r(1-g)$.
\end{itemize}
\end{definition}
\begin{theorem} \label{ampleness of Theta_par} $\cite{pauly1996espaces}$
Let $\eta_*$ a relative family of rank-r parabolic vector bundles of fixed determinant $\delta \in \mathrm{Pic}^{d}(\mathcal{C}/S)$ and of parabolic type $\alpha_*$  parameterized by a $S$-scheme $\T$ over a family of smooth projective curves
$ \mathcal{C}/S$.  Then there is an ample line bundle $\Theta_{par}$  over $ \smp(r,\alpha,\delta)/S$,  such that \quad
$\psi_{\T}^*(\Theta_{par})= \lambda_{par}(\eta_*),$
where $\psi_{\T}$ is the classifying morphism to the relative moduli space $ \smp(r,\alpha,\delta)/S)$.
\end{theorem}

We prove in the following Proposition that, the two definition are the same.
\begin{proposition} \label{B=P} Let $\mathcal{E}_*$ be a family of parabolic rank-r vector bundles of parabolic type $\alpha_*$ over a smooth family of curves $\pi_s: \mathcal{C} \longrightarrow S$ parameterized by a $S$-variety $\T$. Then, \quad $ \Theta_{par}(\eta_*)\cong \lambda_{par}(\eta_*).$
\end{proposition}
\bp To prove the equality of the line bundles over $\T$,  we begin by replacing 
 $\det\left(F^{j}_i(\eta)/F^{j+1}_i(\eta)\right)$ by $ \det\left(\eta_{x_i}/F^{j+1}_i(\eta)\right)$.  In fact we have for all $i \in I$ and $j \in \{1,2,...,\ell_i\}$ the equality
\begin{equation}
\det\left(F^{j}_i(\eta)/F^{j+1}_i(\eta)\right)=\left(\det\left(\eta_{x_i}/F^{j}_i(\eta)\right)\right)^{-1}\otimes \det\left(\eta_{x_i}/F^{j+1}_i(\eta)\right).
\label{transformation}
\end{equation}
for the proof,  we take  for all $i \in I$ and $j \in \{1,2,...,\ell_i\}$,  the quotient  exact sequences 
$$
\xymatrix{ 
0 \ar[r]  & F^{j}_i(\eta) \ar[r] & \eta_{\sigma_i} \ar[r]  &  Q_i^{j-1}(\eta):=  \eta_{\sigma_i}/F^{j}_i(\eta) \ar[r]& 0, \\ 
0 \ar[r]  & F^{j+1}_i(\eta) \ar[r] & \eta_{\sigma_i} \ar[r]  &  Q_i^{j}(\eta):=  \eta_{\sigma_i}/F^{j+1}_i(\eta) \ar[r]& 0, 
}
$$
We calculate the determinants line bundles
\begin{align*}
\det  F^{j}_i(\eta) = \det ( \eta_{\sigma_i}) \otimes \left(\det Q_i^{j-1}(\eta) \right)^{-1} \quad and \quad
 \det F^{j+1}_i(\eta)= \det ( \eta_{\sigma_i}) \otimes \left( \det Q_i^{j}(\eta) \right)^{-1}.
\end{align*}
We calculate determinant using the above equalities,  we get 
\begin{align*}
\det \left(  F^{j}_i(\eta)/  F^{j+1}_i(\eta) \right)&= \left( \det \left(\eta_{\sigma_i}/F^{j}_i(\eta )\right) \right)^{-1} \otimes \det\left(\eta_{\sigma_i}/F^{j+1}_i (\eta)\right).
\end{align*}
Now we can proof the proposition 
\begin{align*}
 \lambda_{par}(\eta_*)& := \lambda(\eta)^{k} \otimes \bigotimes\limits_{i=1}^N \bigotimes_{j=1}^{\ell_i} \left\lbrace \det \left( F^{j}_i(\mathcal{E})/F^{j+1}_i(\mathcal{E})\right) \right\rbrace^{-a_j(i)}  \otimes \det(\mathcal{E}_{\sigma})^{  \frac{k}{r}\chi_{par}} \\
&= \lambda(\eta)^{k} \otimes \bigotimes\limits_{i=1}^N \bigotimes_{j=1}^{\ell_i} \left\lbrace \det \left(\eta_{\sigma_i}/F^{j}_i(\eta )\right)^{-1} \otimes \det\left(\eta_{\sigma_i}/F^{j+1}_i (\eta)\right)  \right\rbrace ^{-a_j(i)}\otimes \det(\mathcal{E}_{\sigma})^{  \frac{k}{r} \chi_{par}}
\end{align*}
by rearranging the terms, we get 
\begin{align*}
\lambda_{par}(\eta_*)= \lambda(\eta)^{k} \otimes \bigotimes\limits_{i=1}^N\left\lbrace   \det(\eta_{\sigma_i})^{a_{\ell_i}(i)} \otimes \bigotimes_{j=1}^{\ell_i-1}  \det \left( \eta_{\sigma_i}/F^{j+1}_i(\eta)\right)^{ p_{j}(i)}  \right\rbrace \otimes \det(\mathcal{E}_{\sigma})^{ \frac{k}{r}\chi_{par}} 
\end{align*}
as $ \det(\mathcal{E}_{\sigma})$ is independent of the section $\sigma$,  we get 
\begin{align*}
\lambda_{par}(\eta_*)= \lambda(\eta)^{k} \otimes \bigotimes\limits_{i=1}^N \bigotimes_{j=1}^{\ell_i-1}  \det \left( \eta_{\sigma_i}/F^{j+1}_i(\eta)\right)^{ p_{j}(i)} \otimes \det(\mathcal{E}_{\sigma})^{\left( \frac{k}{r} \chi_{par}-\sum\limits_{i=1}^N a_{\ell_i}(i)\right)}
\end{align*}

\noindent Now we observe the following equality:
\begin{equation}
\sum\limits_{i=1}^N \sum\limits_{j=1}^{\ell_i} a_j(i)m_j(i)=-\sum\limits_{i=1}^N \sum\limits _{j=1}^{\ell_i-1} d_j(i) r_j(i)+r \sum\limits_{i=1}^N  a_{\ell_i}(i)
 \label{equality}
 \end{equation}
 So the exponent
\begin{align*}
\frac{k}{r} \chi_{par}-\sum\limits_{i=1}^N a_{\ell_i}(i)&=\frac{k}{r} \chi+ \frac{1}{r} \left( \sum\limits_{i=1}^N \sum\limits_{j=1}^{\ell_i}  a_j(i) m_j(i)-r\sum\limits_{i=1}^N a_{\ell_i}(i)  \right)=\frac{k }{r}\chi- \frac{1}{r}\sum\limits_{i=1}^N \sum\limits _{j=1}^{\ell_i-1} d_j(i) r_j(i).
%
\end{align*}
which is equivalent to the equality $k \chi_{par}-r \sum\limits_{i=1}^N  a_{\ell_i}(i)=re.$  This conclude the proof.  \ep \\

We give an other description of the parabolic determinant line bundle.
\begin{proposition} [Parabolic determinant bundle and Hecke modifications] \label{Parabolic determinant bundle and Hecke modifications}
Let $\mathcal{E}_*$ be a family of parabolic rank-r vector bundles of parabolic type $\alpha_*$ over a smooth family of curves $ \mathcal{C}/S$ parameterized by a $S$-variety $\T$,  then 
\begin{equation} 
  \lambda_{par}(\mathcal{E}_*)^r=\Theta^{a} \otimes \left( \bigotimes\limits_{i=1}^N\bigotimes\limits_{j=1}^{\ell_i-1} \Theta_{j}(i)^{n_j(i)p_j(i)} \right),
  \label{decomposition}
 \end{equation} 
where,  for all $i \in I$ and $j \in \{1,2,...,\ell_i-1 \}$
\begin{itemize}
\item $\Theta$ is the pullback of the ample generator of $\mathrm{Pic}(\mathcal{SU}_{\mathcal{C}}(r,\delta)/S)$ by the classifying map $\phi_{\T}$  and  $n=\mathrm{gcd}(r,d)$.
\item $\Theta_{j}(i)$ is the pullback of the ample generators of $\mathrm{Pic}(\mathcal{SU}_{\mathcal{C}}(r,\delta_j(i))/S)$ by the classifying maps $\phi^{\T}_{i,j}$  and  $n_j(i)=\mathrm{gcd}(r,d_j(i))$, where $d_j(i)=\deg(\delta_j(i))$.
\item  $p_j(i)=a_{j+1}(i)-a_j(i)$ \ \  and \ \  $a= n \left( k- \sum\limits_{i=1}^N \sum\limits_{j=1}^{\ell_i-1} p_j(i)  \right)$.
\end{itemize}
\end{proposition}
\bp By the Proposition \ref{parabolic_determinant_2}, the proposition is equivalent to the equality $\Theta_{par}(\eta_*)^r=\Theta^{a} \bigotimes\limits_{i=1}^N\bigotimes\limits_{j=1}^{\ell_i-1} \Theta_{j}(i)^{n_j(i)p_j(i)}.$
\noindent By definition we have
$$\Theta_{par}(\eta_*)= \lambda(\eta)^{k} \otimes \left(\bigotimes\limits_{i=1}^N \bigotimes_{j=1}^{\ell_i-1} \left\lbrace \det \left( \eta_{\sigma_i}/F^{j+1}_i(\eta)\right)^{ p_{j}(i)} \right\rbrace \right) \otimes \det(\mathcal{E}_{\sigma})^{e}$$
 take the Hecke exact sequences 
$$
\xymatrix{ 
0 \ar[r]  & G^{j}_i(\eta) \ar[r] & \eta \ar[r]  &  Q_i^{j}(\eta):=\eta_{\sigma_i}/F^{j+1}_i(\eta) \ar[r]& 0, 
}
$$
by Lemma 3.3 \cite{pauly1998fibres},  we get $\lambda(\eta)=\lambda( G^{j}_i(\eta)) \otimes \left( \det \left(\eta_{\sigma_i}/F^{j+1}_i(\eta) \right) \right)^{-1}.$
\\

We rearrange the terms and by the above equality and take the $r$-th power
\begin{align*}
\Theta_{par}(\eta_*)^r&= \lambda(\eta)^{rk}  \otimes \left( \bigotimes\limits_{i=1}^N \bigotimes_{j=1}^{\ell_i-1} \left\lbrace  \lambda \left( G^{j}_i (\eta) \right) \otimes \lambda(\eta)^{-1} \right\rbrace^{rp_{j}(i)}  \right)\otimes \det(\eta_{\sigma})^{re} \\
&=  \left\lbrace \lambda(\eta)^{k}  \bigotimes\limits_{i=1}^N \bigotimes\limits_{j=1}^{\ell_i-1} \lambda(\eta)^{-p_j(i)} \right\rbrace^r \otimes \left( \bigotimes\limits_{i=1}^N \bigotimes_{j=1}^{\ell_i-1}   \lambda \left( G^{j}_i ( \eta) \right)^{rp_{j}(i)}\right) \otimes \det(\eta_{\sigma})^{re} \\
&= \left\lbrace  \lambda(\eta)^{\frac{r}{n}} \otimes  \det(\mathcal{E}_{\sigma})^{\aleph} \right\rbrace^{a}  \otimes \left( \bigotimes\limits_{i=1}^N \bigotimes_{j=1}^{\ell_i-1}  \left\lbrace  \lambda( G^{j}_i \left( \eta)\right)^{\frac{r}{n_j(i)}} \otimes \det( \eta_{\sigma})^{\aleph_j(i)}  \right\rbrace^{n_j(i) p_{j}(i)}\right)  \otimes \det(\eta_{\sigma})^q.
\end{align*}
where $n \aleph=\chi:=d+r(1-g)$ for $n=\gcd(r,d)$, and  $n_j(i) \aleph_j(i)=\chi_j(i):=d_j(i)+r(1-g)$ for $n_j(i)=\gcd(r,d_j(i))$.
$$\Theta_{par}(\eta_*)^r= \Theta^{a} \otimes  \bigotimes\limits_{i=1}^N \bigotimes_{j=1}^{\ell_i-1} \Theta_j(i)^{n_j(i) p_{j}(i)}\otimes \det(\eta_{\sigma})^q.$$
where
  $a =n \left( k- \sum\limits_{i=1}^N \sum\limits_{j=1}^{\ell_i-1} p_j(i)  \right) \ \ \mathrm{and} \ \
  re= k\chi -\sum\limits_{i=1}^N \sum\limits _{j=1}^{\ell_i-1} p_j(i) r_j(i). $ \\
Hence,  we get $q =re- \aleph a- \left( \sum\limits_{i=1}^N \sum\limits_{j=1}^{\ell_i-1}   \aleph_j(i) p_j(i) n_j(i)\right)=0  $.  This conclude the proof.
\ep 

\paragraph{Canonical bundle} We calculate the canonical bundle of the relative moduli space of rank-r semi-stable parabolic bundles for a fixed parabolic type $\smp_{\mathcal{C}}:=\smp_{\mathcal{C}/S}(r,\alpha_*,\delta)/S $ over a smooth family of marked projective curves parameterized by a scheme $S$.
\paragraph{Canonical bundle in the Grassmannian case:} We suppose that the divisor is of degree one and the flag type is of length one.  Let $\mathcal{E}_*$ a virtual universal parabolic vector bundle 
$$ \def\cartesien{%
    \ar@{-}[]+R+<6pt,-1pt>;[]+RD+<6pt,-6pt>%
    \ar@{-}[]+D+<1pt,-6pt>;[]+RD+<6pt,-6pt>%
  }
\xymatrix{ 
\mathcal{X}^{par} \ar[rr]^{\pi_{n}} \ar[d]_{\pi_{w}} \cartesien &&  \smp_{\mathcal{C}} \ar[d]^{\pi_{e}}  \ar[rr]^{\phi}  && \mathcal{SU}_{\mathcal{C}}(r,\delta) \ar[dll]^{p_e}    \\ 
\ \left( \mathcal{C},D \right) \ar[rr]_{\pi_{s}} && \ar@/^1pc/[ll]^{ \sigma_{i}} S
}
$$
where $D=\sigma(S)$.  In this case the map $\phi$ is a grassmannian bundle over the stable locus of $\mathcal{SU}_{\mathcal{C}}$ (the relative moduli space of semi-stable rank-r vector bundles of determinant $\delta $) and we set $\mathcal{D}:=\pi_w^* \left(D \right)=D \times_S \smp_{\mathcal{C}}$ then we have the Hecke  exact sequence
\begin{equation}
0\longrightarrow \H(\mathcal{E}) \longrightarrow \mathcal{E} \longrightarrow Q(\mathcal{E})=\mathcal{E}\vert_{\mathcal{D}}/F(\mathcal{E}) \longrightarrow 0
\label{B}
\end{equation}  
and the natural exact sequence supported over $\mathcal{D}$
\begin{equation} 
0\longrightarrow F(\mathcal{E}) \longrightarrow \mathcal{E}\vert_{\mathcal{D}} \longrightarrow Q(\mathcal{E}) \longrightarrow 0 
\label{A}
\end{equation}
The relative tangent bundle of the fibration $\phi$ is given as follow 
$$T_{\phi}=\mathrm{ Hom} \left( F(\eta),Q(\eta)) \right)=F(\eta)^{-1}\otimes Q(\eta)$$
we set $r':=\r(Q(\eta))=r-\r(F(\eta))$, hence the canonical bundle is given as follow
$$K_{\phi}=\det (T_{\phi})^{-1}= \det(F(\eta))^{r'} \otimes \det( Q(\eta))^{-(r-r')}$$
the short exact sequence \eqref{A} 
 $$\det(\eta_D)=\det(F(\eta))\otimes \det(Q(\eta)),$$
we replace in the previous equation
$$K_{\phi}=\det(\eta_D)^{r'} \otimes  \det(Q(\eta))^{-r}$$
so we apply  lemma 3.3 \cite{pauly1998fibres}, applied to the Hecke modification sequence \eqref{B} gives the equality 
$$\lambda(\eta)=\lambda(\H(\eta))\otimes \det(Q(\eta))^{-1}$$
which implies that 
$$K_{\phi}=\det(\eta_D)^{r'} \otimes \lambda(\H(\eta))^{-r} \otimes \lambda(\eta)^r 
$$
$$
K_{\phi'}=\left[ \lambda(\H(\eta))^{\frac{r}{n'}} \otimes \det(\eta_y)^{\frac{\chi'}{n'}}   \right]^{-n'} \otimes \left[ \lambda(\eta)^{\frac{r}{n}} \times \det(\mathcal{E}_{D})^{\frac{\chi}{n} } \right]^{n} \otimes \det(\eta_D)^{r'-\chi+\chi'}
$$
where:   $n=\mathrm{gcd}(r,\deg(\mathcal{E}))$, $n'=gcd(r,\deg(G_D(\mathcal{E}))$, $\chi=\chi(\mathcal{E})$ and $\chi'=\chi(G_D(\mathcal{E}))$.
$$ \chi'-\chi=\deg(G_D(\mathcal{E}))-\deg(\mathcal{E})=-r' \Longrightarrow r'-\chi+\chi'=0$$
if we note $\Theta_D$ the pull-back of the ample generator of $\mathrm{Pic}(\mathcal{SU}_{\mathcal{C}}(r,\delta')/S)$,   we get:  $K_{\phi'}=\Theta^n \otimes \Theta_D^{-n'},$ hence
\begin{equation}
K_{\smp_{\mathcal{C}}/S}=K_{\mathcal{SU}_{\mathcal{C}}(r,\delta')/S}\otimes K_{\phi'}=\Theta^{-2n} \otimes \left( \Theta^n \otimes \Theta_D^{-n'} \right)=\Theta^{-n} \otimes  \Theta_D^{-n'}.
\label{grassmanian-canonical-bundle}
\end{equation}
\paragraph{General case:} Now we can calculate the relative canonical bundle $K_{\smp_{\mathcal{C}}/S}$ of the moduli space $\smp_{\mathcal{C}/S}(r,\alpha_*,\delta)$. 
\begin{proposition} \label{Canonical bundle}
Let $b=-n \left( 2+\deg(D)- \sum\limits_{i=1}^N \ell_i  \right).$ Then \quad $K_{\smp_{\mathcal{C}}/S}=\Theta^{b} \otimes \left( \bigotimes\limits_{i=1}^N \bigotimes\limits_{j=1}^{\ell_i-1}  \Theta_{j}(i)^{-n_{j}(i)}\right).$
\end{proposition}
\bp Let $\phi:\smp_{\mathcal{C}} \rightarrow \mathcal{SU}_{\mathcal{C}}(r,\delta)$ be the forgetfull map and let denote its relatif canonical bundle by $K_{\phi}$, then we have $K_{\smp_{\mathcal{C}}/S}=K_{\mathcal{SU}_{\mathcal{C}}(r,\delta) /S}\otimes K_{\phi}$ and by Drezet-Narasimhan \ref{Drezet_Narasimhan} we get $K_{\smp_{\mathcal{C}}/S}=\Theta^{-2n}\otimes K_{\phi}$, where $\Theta$ is the pullback of the relative ample generator of $\mathrm{Pic}(\mathcal{SU}_{\mathcal{C}}(r,\delta)/S)$ by $\phi$ and $n=\gcd(r,\deg(\delta))$. Now, as the map $\phi$ is generically a product of a flag varieties  we can decompose the relative canonical bundle $K_{\phi}=\bigotimes\limits_{i=1}^N  K_{\phi(i)}$, where for all $i\in I$ the bundle $K_{\phi(i)}$ is the canonical bundle of a flag variety. Hence, as the flag variety is embedded canonically in a product of Grassmanians and that its canonical bundle is given by the product of the  canonical bundle over the Grassmanians,  then by the equality \eqref{grassmanian-canonical-bundle},  we have  
$$K_{\phi(i)}=  \bigotimes_{j=1}^{\ell_i-1} \left( \Theta^{n} \otimes  \Theta_{j}(i)^{-n_{j}(i)} \right).$$
We replace and rearrange the terms  in the equation $K_{\smp_{\mathcal{C}}/S}=K_{\mathcal{SU}_{\mathcal{C}}(r,\delta) /S}\otimes K_{\phi}$,  we get 
\begin{align*}
K_{\smp_{\mathcal{C}}/S}&= \Theta^{-2n}\otimes \bigotimes\limits_{i=1}^N K_{\phi(i)} =  \Theta^{-2n}\otimes   \bigotimes\limits_{i=1}^N  \bigotimes_{j=1}^{\ell_i-1} \left( \Theta^{n} \otimes  \Theta_{j}(i)^{-n_{j}(i)} \right) 
\end{align*}
This proves the formula. 
\ep
\section{Hitchin connection for parabolic non-abelian theta-functions}
\begin{theorem} \label{VJ-equation} With the same hypothesis.   Take the parabolic symbol map
$\rho_{par}$,  then the parabolic determinant line bundle $\Theta_{par}$ satisfies the van Geemen-de Jong equation. i.e. ,  \quad
$ \mu_{\Theta_{par}} \circ \rho_{par}=- ( k+r) \Phi_{par}$.\end{theorem}
\bp    By Proposition \ref{Welters},  the Theorem is equivalent to the following points
\begin{enumerate}
\item We prove the equality: $\cup [\Theta_{par} ]\circ \rho_{par}=- k\   \Phi^{par}$, so called the metaplectic correction.  By Proposition \ref{Parabolic determinant bundle and Hecke modifications}, Theorem \ref{equalities} and   linearity with respect to the the tensor product,   we get
\begin{align*}
 \cup [ \Theta_{par}^r ] \circ \rho_{par} &=  \cup \left[ \Theta^{a} \otimes \bigotimes\limits_{i=1}^{N}\bigotimes\limits_{j=1}^{\ell_i-1} \Theta_{j}(i)^{n_j(i)p_j(i)} \right] \circ \rho_{par} \\ 
 &=  a  \left( \cup \left[ \Theta \right] \circ \rho_{par} \right)+  \sum\limits_{i=1}^N \sum\limits_{j=1}^{\ell_i-1} n_j(i)  p_j(i)  \left( \cup \left[ \Theta_j(i) \right] \circ \rho_{par} \right)  \\
   & = - \left( \frac{r}{n}a+ \sum\limits_{i=1}^N \sum\limits_{j=1}^{\ell_i-1} n_j(i)  p_j(i)  \left( \frac{r}{n_j(i)} \right)\right)  \Phi_{par}  
   \end{align*}
   thus
\begin{equation}
     \cup [ \Theta_{par} ] \circ \rho_{par}=-\left( \frac{a}{n}+ \sum\limits_{i=1}^N \sum\limits_{j=1}^{\ell_i-1} p_j(i)  \right) \Phi_{par} =-k \ \Phi_{par}
    \label{1-part-VJ}
\end{equation} 
 which follows by  the following identity 
$$\frac{a}{n}+\sum\limits_{i=1}^N \sum\limits_{j=1}^{\ell_i-1} p_j(i) = \left( k-\sum\limits_{i=1}^N \sum\limits_{j=1}^{\ell_i-1} p_j(i)  \right)+\sum\limits_{i=1}^N \sum\limits_{j=1}^{\ell_i-1} p_j(i)=k
$$

\item We prove the equality: $ \cup [K_{\mathcal{M}^{par}/S}] \circ \rho_{par}= 2r \ \Phi_{par}$. 

By Proposition \ref{Canonical bundle}, Theorem \ref{equalities} and linearity with respect to the the tensor product,  we have
\begin{align*}
  \cup [K_{\smp_{\mathcal{C}}/S}] \circ \rho_{par} &= \cup \left[ \Theta^{-b} \otimes \bigotimes\limits_{i=1}^N \bigotimes\limits_{j=1}^{\ell_i-1}  \Theta_j(i)^{-n_j(i)} \right] \circ \rho_{par} \\
  &=-b \left( \cup [ \Theta] \circ \rho_{par} \right) + \sum\limits_{i=1}^N \sum\limits_{j=1}^{\ell_i-1} -n_j(i) \left(  \cup [ \Theta_j(i) ] \circ \rho_{par} \right) \\
  &= r \left( 2+\deg(D)- \sum\limits_{i=1}^N \ell_i+ \sum\limits_{i=1}^N \sum\limits_{j=1}^{\ell_i-1}1  \right)\ \Phi_{par}= 2r \ \Phi_{par}.
\end{align*}
\end{enumerate}
adding the two equations,  we get
 \begin{equation}
 \mu_{\Theta_{par}} \circ \rho_{par} =-(k+r) \cdot \Phi_{par}.
 \label{1+2-part-VJ}
 \end{equation}  
\ep

 We observe that the composition $\mu_{\Theta_{par}}\circ \rho_{par}$ does not depend on the parabolic weights but depends on the level-$k$,  in some sense what contribute in the decomposition \eqref{decomposition} in the term $\Theta^k$,  we rearrange the terms as follow
\begin{align*}
\Theta_{par}^r =\Theta^{a} \bigotimes\limits_{i=1}^N\bigotimes\limits_{j=1}^{\ell_i-1} \Theta_j(i)^{n_j(i)p_j(i)} 
=\Theta^k \otimes \left(  \bigotimes\limits_{i=1}^N\bigotimes\limits_{j=1}^{\ell_i-1} \left( \Theta^{-n} \otimes \Theta_j(i)^{n_j(i)} \right)^{p_j(i)}\right).
\end{align*}
by the Definition \ref{parabolic_determinant_2} and Propositions \ref{B=P}, \ref{Parabolic determinant bundle and Hecke modifications}  the identification we get 
$$ \bigotimes\limits_{i=1}^N \bigotimes_{j=1}^{\ell_i-1} \left\lbrace \det \left( \eta_{\sigma_i}/F^{j+1}_i(\eta)\right)^{ rp_{j}(i)} \right\rbrace \otimes \det(\mathcal{E}_{\sigma})^{(re-k\chi)}=\bigotimes\limits_{i=1}^N\bigotimes\limits_{j=1}^{\ell_i-1} \left( \Theta^{-n} \otimes \Theta_j(i)^{n_j(i)} \right)^{p_j(i)}.$$
which we call the flag part of the determinant line bundle and we denoted $\mathcal{F}(\alpha_*)$.  By Corollary \ref{equalities} for all $i \in I$ and $j \in \{1,2,...,\ell_i \}$   we have 
$$ \left( \cup \left[ \Theta^{-n} \right]+ \cup \left[ \Theta_j(i)^{n_j(i)} \right] \right) \circ \rho_{par} = 0$$
thus by Proposition \ref{rho_isomorphism} we have that $\rho_{par}$ is an isomorphism,  then we have: 
\begin{equation}
\cup \left[ \mathcal{F}(\alpha_*)\right] =0.
\label{flag-equation}
\end{equation}
\begin{remark} 
 \begin{enumerate}
\item In general case (see.  $\cite{singh2021differential}$)  : If $X$ is a Hitchin variety and $L$ line bundle over $X$,  then 
$$\cup \left[L \right]: \mathrm{H}^0(X, \mathrm{Sym}^q \ T_X) \longrightarrow \mathrm{H}^1(X, \mathrm{Sym}^{q-1}T_X)$$
which can be seen as the first connecting map on cohomology of the short exact sequence 
$$ 0 \longrightarrow \mathrm{Sym}^{q-1}\left( \mathcal{D}^{(1)}_X(L)\right) \longrightarrow \mathrm{Sym}^q( \mathcal{D}^{(1)}_X(L)) \longrightarrow \mathrm{Sym}^q(T_X) \longrightarrow 0$$
which is the $q$-th symmetric power of the the Atiyah sequence \eqref{atiyah-classe},  and we have the following theorem
\begin{theorem} [\cite{singh2021differential},  Theorem 2.2] \label{Singh}
If $L$ is an ample line bundle then the map above is an isomorphism.
\end{theorem}
\item The varieties $\mathcal{SU}_{\mathcal{C}/S}(r,\delta)$ and $\smp_{\mathcal{C}/S}(r,\alpha_*,\delta)$ are Hitchin varieties in the sense of $\cite{singh2021differential}$,  and by Theorem \ref{ampleness of Theta_par} the parabolic determinant line bundle $\Theta_{par}$ is ample thus the map is an isomorphism
$$\cup \left[ \Theta_{par} \right]: \pi_{e_*} \mathrm{Sym}^2 \left( T_{\smp/S} \right) \longrightarrow R^1\pi_{e_*} \left( T_{\smp/S}\right),$$
\end{enumerate}
\end{remark}

By equalities  \eqref{1+2-part-VJ} and \eqref{flag-equation},   one has for all positive integer $\nu$ the equalities
\begin{align*}
\mu_{\Theta_{par}^{\nu}}&= \frac{n \left( \nu k+r \right)}{r} \cdot \cup \left[\Theta \right]=\left( \frac{\nu k+r}{k} \right) \cup \left[ \Theta_{par} \right].
\end{align*}

Thus, by the previous remark, we have the following proposition.
 
\begin{proposition} \label{mu-isomorphism} For any positive integer $\nu$, the map $\mu_{\Theta_{par}^{\nu}}$ is an isomorphism.
\end{proposition}

We get the van Geemen-de jong equation Theorem \ref{VJ-equation}  for any positive power of the parabolic theta line bundle 
$$ \mu_{\Theta_{par}^{\nu}} \circ \rho_{par}=- ( \nu k+r) \Phi_{par}  $$
 \begin{theorem} \label{Main Theorem 1}    Consider a smooth versal family $(\pi_s:(\mathcal{C},D) \longrightarrow S)$ of complex projective marked curves of genus $g \geq 2$, and $D$ a reduced divisor of relative degree $N$.  Take $\alpha_*=(k,\vec{a},\vec{m})$ a fixed rank-$r$ parabolic type with respect to the divisor $D$. We denote by $(\pi_e: \smp(r,\alpha_*,\delta)  \rightarrow S)$, the relative moduli space of parabolic rank-$r$ vector bundles over $(\mathcal{C},D)/S$ with determinant $\delta \in \mathrm{Pic}^d(\mathcal{C}/S)$,  equipped with the parabolic determinant line bundle $\Theta_{par}$.Let $\nu$ be a positive integer. Then, there exists a unique projective flat connection on the vector bundle ${\pi_e}_*(\Theta_{par}^{\nu})$ of non-abelian parabolic theta functions,  induced by a heat operator with symbol $$ \rho^{Hit}_{par}(\nu):= \frac{1}{(\nu k+r)} \left( \rho_{par}\circ \kappa_{\mathcal{C}/S} \right).$$  
 \end{theorem}
 
\bp Let prove the theorem for $\nu=1$,  we denote  $ \rho^{Hit}_{par}:= \rho^{Hit}_{par}(1)$.
\begin{itemize}
\item  First we prove existence of the connection:
 We apply van Geemen-de Jong Theorem \ref{van Geemen and De Jong} for $L=\Theta_{par}$ over $\smp_{\mathcal{C}/S}(r,\alpha_*,\delta)$.  Thus by Theorem \ref{VJ-equation} and Proposition \ref{Parabolic BSBE},  we get the first condition  of Theorem  \ref{van Geemen and De Jong}
 \begin{align*}
  \mu_{\Theta_{par}} \circ \left( \rho_{par} \circ \kappa_{\mathcal{C}/S} \right) &= -(k+r) \Phi_{par} \circ \kappa_{\mathcal{C}/S} =-(k+r)\ \kappa_{\smp_{\mathcal{C}}/S}.
  \end{align*}
  The second condition follows by Theorem \ref{Singh} for $q=1$.i.e.  the map
 $$\cup \left[ \Theta_{par} \right]: \pi_{e_*}  (T_{\smp/S}) \longrightarrow R^1\pi_{e_*}( \mathcal{O}_{\smp/S}) $$
 is an isomorphism, as the relative Picard group $\mathrm{Pic}(\smp_{\mathcal{C}}/S)$ is discrete then the infinitesimal deformations of any line bundle $L$ over $\smp$ are trivial and parameterized by the sheaf  $\mathrm{R}^1\pi_{e_*}(\mathcal{O}_{\smp})$ over $S$,  thus  $\pi_{e_*}( \mathcal{O}_{\smp})\cong 0$, as consequence $\pi_{e_*} T_{\smp/S} \cong 0$, no global relative vector fields. 
  
The third condition follows from the algebraic Hartogs's Theorem and the fact that the space $\smp_{\mathcal{C}}(r,\alpha_*,\delta)$ is normal  variety and proper over $S$ and that the smooth locus is a	 big open subset of  $\smp_{\mathcal{C}}(r,\alpha_*,\delta)$.
  \item Flatness of the connection: We apply flatness criterion Theorem  \ref{Flatness criterion} to the parabolic symbol map $\rho_{par}$. 
The first condition holds since by definition of the parabolic Hitchin map corresponds to homogeneous functions on $T^{\vee}_{\smp/S}$ of degree two under the action of $\mathbb{G}_m$ in the parabolic Hitchin system,  hence Poisson-commute.  The second point is given in Proposition \ref{mu-isomorphism} and the third point is is given in the first part of the proof.
  \end{itemize} 
  \ep
  
For $D=\emptyset$ and $\alpha_*=k \in \NN^*$ the trivial parabolic type.   We have the identification $\smp_{\mathcal{C}}(r,0_*,\delta)=\mathcal{SU}_{\mathcal{C}}(r,\delta)$ the moduli space of semistable rank-r vector bundles with fixed determinant $\delta \in \mathrm{Pic}^d(\mathcal{C}/S)$,  hence $\Theta_{par}^{r/n}=\l^k$ for $n:=\gcd(r,\deg(\delta))$ and $\rho_{par}\equiv \rho^{Hit}$.  Then as a corollary we obtain the following theorem
\begin{theorem} \label{Main Theorem 2}
 We consider $\mathcal{C}/S$  a smooth family  of complex projective curves of genus $g \geq 2$ (and $g \geq 3$ if $r=2$ and $\deg(\delta)$ is even). Let $\l$ be the determinant line bundle over $p_e: \mathcal{SU}_{\mathcal{C}/S}(r,\delta) \rightarrow S$. Let $k$ be a positive integer. Then, there exists a unique projective flat connection on the vector bundle $p_{e_*}(\l^{k})$ of non-abelian theta functions,  induced by a heat operator with symbol map
$$ \rho(k):= \frac{n}{r(k+n)} \left( \rho^{Hit}\circ \kappa_{\mathcal{C}/S} \right).$$
 \end{theorem}
 \paragraph{Some comments:}
We apply Theorem \ref{Main Theorem 2} for $\delta=\mathcal{O}_{\mathcal{C}}$ thus $n=r$ and $\Theta_{par}=\l^k$,  we get 
$$ \rho(k):= \frac{1}{(k+r)} \left( \rho_{Hit}\circ \kappa_{\mathcal{C}/S} \right)$$
Hence, we recover Theorems  4.8.1 and 4.8.2 in \cite{pauly2023hitchin},  which was generalized in \cite{biswas2021ginzburg,biswas2021geometrization}  to $\smp_{\mathcal{C}/S}(r,\alpha_*)$ the space of rank-r parabolic bundle with trivial determinant and parabolic type $\alpha_*$.  The symbol map is given for a positive integer $\nu$ by 
$$\rho^{Hit}_{par, \Gamma}(\nu):=\vert \Gamma \vert \ \mu^{-1}_{\Theta^{\nu}} \circ \left(\cup[\Theta] \circ \rho_{par} \circ \kappa_{\mathcal{C}/S} \right),$$
By \ref{mu-isomorphism} and  \eqref{flag-equation},  we get the equality:  $\vert \Gamma\vert (\mu^{-1}_{\Theta_{par}^{\nu}} \circ \cup[\Theta]) = \frac{\vert \Gamma \vert}{(\nu k+r)} \mathrm{Id}$,  hence
 $$\rho^{Hit}_{par, \Gamma}(\nu):= \frac{\vert \Gamma \vert}{(\nu k+r)}    \left(  \rho_{par} \circ \kappa_{\mathcal{C}/S} \right)=\vert \Gamma \vert \rho^{Hit}_{par}(\nu).$$
 The factor $\vert \Gamma \vert$ is because they work over $\mathcal{SU}_{X/S}^{\Gamma}(r)$ the moduli space of $\Gamma$-linearised bundles for family of Galois coverings $h:X \longrightarrow (\mathcal{C},D)$ parameterized by the variety $S$. 
\begin{remark} If the system of weights $\alpha_*$ is not generic in the sense of Yokogawa,  then the moduli space $\mathcal{SM}^{par}_{\mathcal{C}/S}(r,\alpha_*,\delta)$ is not smooth,  and its Picard group is not maximal.  In other words,  not all line bundles on the Quot scheme descend to the moduli space.   In fact,  we can choose the weights $\alpha_*$ in such a way that we have the following isomorphism: $$\mathrm{Pic}\left( \mathcal{SM}^{par}_{\mathcal{C}/S}(r,\alpha_*,\delta)/S \right) \simeq \mathbb{Z} \Theta_{par}.$$

To prove Theorem \ref{Main Theorem 1} in the case of non-generic weights, we must work over the stack of quasi-parabolic vector bundles,  where the Picard group is maximal,  the Hecke maps $\H_i^j$ and the forgetful map are maps between stacks (no stability conditions).  Note that the decompositions of the parabolic determinant line bundle \ref{Parabolic determinant bundle and Hecke modifications} and the canonical line bundle \ref{Canonical bundle} still hold.  Hence, we get the van Geemen-de Jong equation over the stack,  and as the morphisms that appear in this equation are well defined over the moduli space of parabolic bundles for any system of weights,  we get the existence theorem of the Hitchin connection. 
\begin{example}[Non-generic weights.  See \cite{pauly1998fibres} for the details] Let's consider the rank two case with $D$ being a parabolic divisor of degree $N=2m \geq 4$.  For  all $i\in I$,  we choose the following system of weights:
\begin{center}
$a_1(i)=0 $, $a_2(i)=1$ and $k=2$.    
\end{center}
In this case, the Picard group of the moduli space $\mathcal{SM}^{par}_{\mathcal{C}/S}(2,\alpha_*,\mathcal{O}_{\mathcal{C}})$ is generated by the line bundle $\Theta_{par}$.
\end{example}
\end{remark}
 
\bibliographystyle{alpha}

\begin{thebibliography}{BMW21b}

\bibitem[Ati57]{atiyah1957complex}
M.  Atiyah.
\newblock Complex analytic connections in fibre bundles.
\newblock {\em Transactions of the American Mathematical Society},
  85(1):181--207, 1957.
  
\bibitem[AB83]{Atiyah-Bott}
M.Atiyah and R.~Bott.
\newblock The Yang-Mills equations over Riemann surfaces.
\newblock {\em Philos. Trans. Roy. Soc. London Ser. A}, 308(1505):523--615,
  1983.

\bibitem[AGL12]{andersen2012hitchin}
J. Andersen,  N.  Gammelgaard,  and M.
  Lauridsen.
\newblock Hitchin's connection in metaplectic quantization.
\newblock {\em Quantum topology}, 3(3):327--357, 2012.

\bibitem[And12]{andersen2012hitchinT}
J. Andersen.
\newblock Hitchin's connection,  Toeplitz operators, and symmetry invariant deformation quantization.
\newblock {\em Quantum topology}, 3(3):293--325, 2012.

\bibitem[Axe91]{axelrod1991geometric}
S.  Axelrod.
\newblock {\em Geometric quantization of Chern-Simons gauge theory}.
\newblock Princeton University, 1991.

\bibitem[BBMP23]{pauly2023hitchin}
T.  Baier,  M.  Bolognesi,  J. Martens, and C.  Pauly.
\newblock The Hitchin connection in arbitrary characteristic.
\newblock {\em Journal of the Institute of Mathematics of Jussieu},
  22(1):449--492, 2023.

\bibitem[BBN03]{MR1995862}
V.~Balaji, I.~Biswas, and D.~S. Nagaraj.
\newblock Ramified {$G$}-bundles as parabolic bundles.
\newblock {\em J. Ramanujan Math. Soc.}, 18(2):123--138, 2003.

\bibitem[Bau89]{bauer1989parabolic}
S.~Bauer.
\newblock Parabolic bundles, elliptic surfaces and SU(2)-representation spaces
  of genus zero Fuchsian groups.
\newblock 1989.

\bibitem[BL94]{beauville1994conformal}
A.  Beauville and Y.  Laszlo.
\newblock Conformal blocks and generalized theta functions.
\newblock {\em Communications in mathematical physics}, 164:385--419, 1994.

\bibitem[BS88]{beilinson1988determinant}
A.  Beilinson and V.  Schechtman.
\newblock Determinant bundles and virasoro algebras.
\newblock {\em Communications in mathematical physics}, 118:651--701, 1988.

\bibitem[Ber93]{bertram1993generalized}
A.  Bertram.
\newblock Generalized  SU(2)-theta functions.
\newblock {\em Inventiones mathematicae}, 113(1):351--372, 1993. 

\bibitem[Bho89]{bhosle1989parabolic}
U.  Bhosle.
\newblock Parabolic vector bundles on curves.
\newblock {\em Arkiv f{\"o}r matematik}, 27(1-2):15--22, 1989.

\bibitem[BE99]{bloch1999relative}
S.  Bloch and H.  Esnault.
\newblock Relative algebraic differential characters.
\newblock {\em arXiv preprint math/9912015}, 1999.

\bibitem[Biq91]{biquard1991fibres}
O. Biquard.
\newblock Fibr{\'e}s paraboliques stables et connexions singulieres plates.
\newblock {\em Bulletin de la Soci{\'e}t{\'e} Math{\'e}matique de France},
  119(2):231--257, 1991.

\bibitem[BMW21]{biswas2021geometrization}
I.  Biswas,  S. Mukhopadhyay, and R. Wentworth.
\newblock Geometrization of the tuy/wzw/kz connection.
\newblock {\em arXiv preprint arXiv:2110.00430}, 2021.

\bibitem[BDHP22]{biswas2022infinitesimal}
I.  Biswas,  S.  Dumitrescu,  S. Heller,  and C.  Pauly.
\newblock Infinitesimal deformations of parabolic connections and parabolic opers.
\newblock {\em arXiv preprint arXiv:2202.09125}, 2022.

\bibitem[BMW23]{biswas2021ginzburg}
I.  Biswas,  S. Mukhopadhyay, and R. Wentworth.
\newblock A Hitchin connection on nonabelian theta functions for parabolic G-bundles.
\newblock {\em Journal f\"{u}r die reine und angewandte Mathematik (Crelles Journal)} (0), 2023.

\bibitem[BMW23b]{biswas2023BS}
I.  Biswas,  S. Mukhopadhyay, and R. Wentworth.
\newblock A parabolic analog of a theorem of Beilinson and Schechtman.
\newblock {\em  arXiv preprint arXiv:2307.09196},  2023.

\bibitem[BR93]{BR93}
I.~Biswas and N.~Raghavendra.
\newblock Determinants of parabolic bundles on {R}iemann surfaces.
\newblock {\em Proc. Indian Acad. Sci. Math. Sci.}, 103(1):41--71, 1993.

\bibitem[BY99]{boden1999rationality}
H. Boden and K.  Yokogawa.
\newblock Rationality of moduli spaces of parabolic bundles.
\newblock {\em Journal of the London Mathematical Society}, 59(2):461--478,
  1999.

\bibitem[B18]{bjerre2018hitchin}
M.  Bjerre.
\newblock The Hitchin connection for the quantization of the moduli space of  parabolic bundles on surfaces with marked points.
\newblock {\em PhD Dissertations, Department of Mathematics, Aarhus
  University}, 2018.

\bibitem[CS07]{crampin2007projective}
M.  Crampin and DJ~Saunders.
\newblock Projective connections.
\newblock {\em Journal of Geometry and Physics}, 57(2):691--727, 2007.

\bibitem[DN89]{drezet1989groupe}
J.M. Drezet and M.S. Narasimhan.
\newblock Groupe de picard des vari{\'e}t{\'e}s de modules de fibr{\'e}s
  semi-stables sur les courbes alg{\'e}briques.
\newblock {\em Inventiones mathematicae}, 97:53--94, 1989.

\bibitem[ET00]{esnault2000determinant}
H.~Esnault and I-Hsun Tsai.
\newblock Determinant bundle in a family of curves after Beilinson and Schechtman.
\newblock {\em Communications in Mathematical Physics}, 211(2):359--363, 2000.

\bibitem[Fal93]{faltings1993stable}
G.~Faltings.
\newblock Stable G-bundles and projective connections.
\newblock {\em J. Algebraic Geom}, 2(3):507--568, 1993.

\bibitem[FGI{\etalchar{+}}05]{fantechi2005fundamental}
B.  Fantechi,  L.  G{\"o}ttsche,  L.  Illusie,  S. Kleiman,  N.
  Nitsure, and A.  Vistoli.
\newblock {\em Fundamental algebraic geometry: Grothendieck's FGA explained},
  volume 123.
\newblock American Mathematical Society Providence, RI, 2005.

\bibitem[GdJ98]{van1998hitchin}
B.~Van Geemen and A.J de~Jong.
\newblock On Hitchin's connection.
\newblock {\em Journal of the American Mathematical Society}, 11(1):189--228,
  1998.

\bibitem[Har13]{hartshorne2013algebraic}
R.~Hartshorne.
\newblock {\em Algebraic geometry}, volume~52.
\newblock Springer Science \& Business Media, 2013.

\bibitem[Hit87]{hitchin1987stable}
N~Hitchin.
\newblock Stable bundles and integrable systems.
\newblock {\em Duke mathematical journal}, 54(1):91--114, 1987.

\bibitem[Hit90a]{hitchin1990flat}
N. Hitchin.
\newblock Flat connections and geometric quantization.
\newblock {\em Communications in mathematical physics}, 131:347--380, 1990.

\bibitem[Hit90b]{hitchin1990symplectic}
N.J. Hitchin.
\newblock The symplectic geometry of moduli spaces of connections and geometric
  quantization.
\newblock {\em Progress of Theoretical Physics Supplement}, 102:159--174, 1990.

\bibitem[KD76]{knudsen1976projectivity}
F.~Knudsen and D.Mumford.
\newblock The projectivity of the moduli space of stable curves i:
  Preliminaries on" det" and" div".
\newblock {\em Mathematica Scandinavica}, 39(1):19--55, 1976.

\bibitem[Lan61]{lang1961grothendieck}
S~Lang.
\newblock A.  Grothendieck, {\'e}l{\'e}ments de g{\'e}om{\'e}trie
  alg{\'e}brique.
\newblock {\em Bull. Amer. Math. Soc.}, 67(6):239--246, 1961.

\bibitem[Las98]{laszlo1998hitchin}
Y.  Laszlo.
\newblock Hitchin's and wzw connections are the same.
\newblock {\em Journal of Differential Geometry}, 49(3):547--576, 1998.

\bibitem[LS97]{laszlo1997line}
Y.~Laszlo and C.~Sorger.
\newblock The line bundles on the moduli of parabolic $G$-bundles over curves and their sections.
\newblock In {\em Annales scientifiques de l'Ecole normale sup{\'e}rieure},
  volume~30, pages 499--525, 1997.

\bibitem[Mar09]{martinengo2009higher}
E.~Martinengo.
\newblock {\em Higher brackets and moduli space of vector bundles}.
\newblock PhD thesis,  Sapienza Universita di Roma, 2009.

\bibitem[MY92]{maruyama1992moduli}
M.~Maruyama and K.~Yokogawa.
\newblock Moduli of parabolic stable sheaves.
\newblock {\em Mathematische Annalen}, 293:77--99, 1992.

\bibitem[MS80]{mehta1980moduli}
V. Mehta and C. Seshadri.
\newblock Moduli of vector bundles on curves with parabolic structures.
\newblock {\em Math. Ann}, 248:205--239, 1980.

\bibitem[NR93]{narasimhan1993factorisation}
M.  Narasimhan and T.  Ramadas.
\newblock Factorisation of generalised theta functions. i.
\newblock {\em Inventiones mathematicae}, 114(1):565--623, 1993.

\bibitem[Pau96]{pauly1996espaces}
C.  Pauly.
\newblock Espaces de modules de fibr{\'e}s paraboliques et blocs conformes.
\newblock {\em Duke Math. J.}, 85(1):217--235, 1996.

\bibitem[Pau98]{pauly1998fibres}
C.Pauly.
\newblock Fibr{\'e}s paraboliques de rang 2 et fonctions th{\^e}ta
  g{\'e}n{\'e}ralis{\'e}es.
\newblock {\em Mathematische Zeitschrift}, 228:31--50, 1998.

\bibitem[Qui85]{MR783704}
D.~Quillen.
\newblock Determinants of {C}auchy-{R}iemann operators on {R}iemann surfaces.
\newblock {\em Funktsional. Anal. i Prilozhen.}, 19(1):37--41, 96, 1985.

\bibitem[Ran06]{ran2006jacobi}
Z. Ran.
\newblock Jacobi cohomology, local geometry of moduli spaces,  and Hitchin
  connections.
\newblock {\em Proceedings of the London Mathematical Society}, 92(3):545--580,
  2006.

\bibitem[Sch13]{schaposnik2013spectral}
L.  Schaposnik.
\newblock Spectral data for G-Higgs bundles.
\newblock {\em arXiv preprint arXiv:1301.1981}, 2013.

\bibitem[SS95]{schottenloher1995metaplectic}
P.  Scheinost and M.  Schottenloher.
\newblock Metaplectic quantization of the moduli spaces of flat and parabolic
  bundles.
\newblock {\em j. reine angew. Math}, 466:145--219, 1995.

\bibitem[Ser07]{sernesi2007deformations}
E.~Sernesi.
\newblock {\em Deformations of algebraic schemes}, volume 334.
\newblock Springer Science \& Business Media, 2007.

\bibitem[Ses77]{seshadri1977moduli}
C.~Seshadri.
\newblock Moduli of vector bundles with parabolic structures.
\newblock {\em Bull. Amer. Math. Soc.}, 83, 1977.

\bibitem[Ses82]{seshadri1982fibers}
C.~Seshadri.
\newblock Fib\'ers vectoriels sur les courbes alg\'ebriques.
\newblock {\em Ast{\'e}risque}, 96:1--209, 1982.

\bibitem[Ses11]{seshadri2011moduli}
C.~Seshadri.
\newblock Moduli of $\pi$-vector bundles over an algebraic curve.
\newblock {\em Questions on algebraic varieties}, pages 139--260, 2011.

\bibitem[Sim90]{simpson1990harmonic}
C.~Simpson.
\newblock Harmonic bundles on noncompact curves.
\newblock {\em Journal of the American Mathematical Society}, 3(3):713--770,
  1990.

\bibitem[Sin21]{singh2021differential}
A.~Singh.
\newblock Differential operators on Hitchin variety.
\newblock {\em Journal of Algebra}, 566:361--373, 2021.

\bibitem[ST04]{sun2004hitchin}
X.~Sun and I.~Tsai.
\newblock Hitchin's connection and differential operators with values in the  determinant bundle.
\newblock {\em Journal of Differential Geometry}, 67(2):335--376, 2004.

\bibitem[TUY89]{tsuchiya1989conformal}
A.~Tsuchiya,  K.~Ueno,  and Y.~Yamada.
\newblock Conformal field theory on universal family of stable curves with gauge symmetries.
\newblock In {\em Integrable Sys Quantum Field Theory}, pages 459--566.
  Elsevier, 1989.

\bibitem[Wel83]{welters1983polarized}
G.E Welters.
\newblock Polarized abelian varieties and the heat equations.
\newblock {\em Compositio Mathematica}, 49(2):173--194, 1983.

\bibitem[Wit89]{witten1989quantum}
E.~Witten.
\newblock Quantum field theory and the jones polynomial.
\newblock {\em Communications in Mathematical Physics}, 121(3):351--399, 1989.

  \bibitem[Yok91]{yokogawa1991moduli}
K.~Yokogawa.
\newblock Moduli of stable pairs.
\newblock {\em Journal of Mathematics of Kyoto University}, 31(1):311--327,
  1991.
  
\bibitem[Yok93]{yokogawa1993compactification}
K.~Yokogawa.
\newblock Compactification of moduli of parabolic sheaves and moduli of parabolic Higgs sheaves.
\newblock {\em Journal of Mathematics of Kyoto University}, 33(2):451--504,
  1993.

\bibitem[Yok95]{yokogawa1995infinitesimal}
K.~Yokogawa.
\newblock Infinitesimal deformation of parabolic higgs sheaves.
\newblock {\em International Journal of Mathematics}, 6(1):125, 1995.
\end{thebibliography}
\nocite{*}
{\footnotesize
\newcommand{\etalchar}[1]{$^{#1}$}

}
\end{document}